\documentclass[a4paper]{scrartcl}
\usepackage[utf8]{inputenc}
\usepackage{natbib}
\usepackage{amsthm, amssymb, amsmath}
\usepackage{bm}
\usepackage{dsfont}
\usepackage{enumerate}
\usepackage{color}
\usepackage{soul}
\usepackage{caption}
\usepackage{graphicx}
\usepackage[FIGTOPCAP]{subfigure}
\usepackage[unicode,breaklinks]{hyperref}

\definecolor{linkcolor}{RGB}{15,15,175}
\definecolor{alarmcolor}{RGB}{175,15,15}

\hypersetup{
    pdfauthor={Christian Sch\"afer},
    pdfcreator={Christian Sch\"afer},
    unicode=true,           
    pdftoolbar=true,        
    pdfmenubar=true,        
    pdffitwindow=false,     
    pdfstartview={FitH},    
    pdfnewwindow=true,      
    colorlinks=true,        
    linkcolor=linkcolor,    
    citecolor=linkcolor,    
    filecolor=linkcolor,    
    urlcolor =linkcolor     
}

\def\B{\mathbb{B}}

\def\R{\mathbb{R}}

\def\ber{\mathcal{B}} 
\def\nor{\mathcal{N}} 

\def\diag#1{\mathrm{diag}\left[ #1 \right]}

\def\abs#1{\left\vert #1 \right\vert}

\def\t{^\intercal}

\def\v#1{\bm{#1}}
\def\m#1{\bm{\mathrm{#1}}}

\def\eqdef{:=}
\def\proptodef{:\propto}

\def\ind{\mathds{1}}
\def\logit{\ell}

\newcommand{\normal}[2]{\nor\left(#1,#2\right)}
\newcommand{\prob}[1]{\mathbb{P}\left(#1 \right)}
\newcommand{\probplus}[2]{\mathbb{P}_{#1}\left(#2 \right)}

\newcommand{\evplus}[2]{\mathbb{E}_{#1}\left(#2 \right)}

\def\keywords#1{\par\addvspace\medskipamount{\rightskip=0pt plus1cm
\def\and{\ifhmode\unskip\nobreak\fi\ $\cdot$
}\noindent \textbf{Keywords}\enspace\ignorespaces#1\par}}

\addtokomafont{footnote}{\sffamily\fontsize{6}{3}\selectfont}
\author{Christian Sch\"afer$^{1,2}$}
\def\crest{
\footnotetext[1]{Centre de Recherche en Économie et Statistique, 3 Avenue Pierre Larousse, 92240 Malakoff, France}
\footnotetext[2]{CEntre de REcherches en MAthématiques de la DEcision, Université Paris-Dauphine, Place du Maréchal de Lattre de Tassigny
75775 Paris, France}
}

\usepackage{setspace}
\usepackage{algorithmic}
\usepackage{algorithm}

\floatname{algorithm}{Procedure}

\def\up#1{^{(#1)}}

\author{Christian Sch\"afer$^{1,2}$\hspace{1cm} Nicolas Chopin$^{1,3}$}
\title{Sequential Monte Carlo on large binary sampling spaces}

\hypersetup{
    pdfauthor={Christian Sch\"afer, Nicolas Chopin},
    pdfcreator={Christian Sch\"afer},
    pdftitle={Sequential Monte Carlo on large binary sampling spaces},
    pdfkeywords={Adaptive Monte Carlo} {Multivariate binary data} {Sequential Monte Carlo} {Bayesian Variable selection}
}

\begin{document}
\maketitle
\crest
\thispagestyle{empty}
\footnotetext[3]{Ecole Nationale de la Statistique et de l'Administration, 3 Avenue Pierre Larousse, 92240 Malakoff, France}

%

%


%

%

\begin{abstract}
A Monte Carlo algorithm is said to be adaptive if it automatically calibrates its current proposal distribution using past simulations. The choice of the parametric family that defines the set of proposal distributions is critical for good performance. In this paper, we present such a parametric family for adaptive sampling on high-dimensional binary spaces.

A practical motivation for this problem is variable selection in a linear regression context. We want to sample from a Bayesian posterior distribution on the model space using an appropriate version of Sequential Monte Carlo.

Raw versions of Sequential Monte Carlo are easily implemented using binary vectors with independent components. For high-dimensional problems, however, these simple proposals do not yield satisfactory results. The key to an efficient adaptive algorithm are binary parametric families which take correlations into account, analogously to the multivariate normal distribution on continuous spaces.

We provide a review of models for binary data and make one of them work in the context of Sequential Monte Carlo sampling. Computational studies on real life data with about a hundred covariates suggest that, on difficult instances, our Sequential Monte Carlo approach clearly outperforms standard techniques based on Markov chain exploration.
\end{abstract}

\keywords{
Adaptive Monte Carlo \and Multivariate binary data \and Sequential Monte Carlo \and Linear regression \and Variable selection
}

%

%


%

%

\thispagestyle{empty}
\section{Introduction}
\label{sec:intro}
We present a Sequential Monte Carlo \citep{del2006sequential} algorithm for adaptive sampling from a binary distribution. A Monte Carlo algorithm is said to be adaptive if it adjusts, sequentially and automatically, its sampling distribution to the problem at hand. Besides Sequential Monte Carlo, important classes of adaptive Monte Carlo are Adaptive Importance Sampling \citep[e.g.][]{cappe2008adaptive} and Adaptive Mar\-kov chain Monte Carlo \citep[e.g.][]{andrieu2008tutorial}.

A central aspect of adaptive algorithms is their need for a parametric family of auxiliary distributions which should have the following three properties:
(a)
the family is sufficiently flexible to guarantee a reasonable performance in the context of the specific algorithm;
(b)
it allows to quickly draw independent samples;
(c)
it can, with reasonable effort, be calibrated using past simulations.

For problems in continuous sampling spaces, the typical example is the multivariate normal distribution, which clearly fulfils (b) and (c), and complies with (a) in many practical problems. In this paper, we propose an analogue for high-dimensional binary sampling spaces. 

\subsection{Adaptive Monte Carlo on multivariate binary spaces}
\label{sec:amc ob bs}
Our objective is to construct a parametric family for Sequential Monte Carlo on the binary sampling space $\B^d=\lbrace0,1\rbrace^d$, where $d$ is too large to allow for exhaustive enumeration of the whole space $\B^d$. Since there is no multivariate binary family which we can easily parametrise by its first and second order moments like the multivariate normal, the construction of suitable proposal distributions seems more difficult for the discrete adaptive sampling problem than for its continuous counterpart.

The major application for our algorithm is variable selection in linear regression models. In this context, a binary vector $\v \gamma\in\B^d$ encodes whether each of $d$ possible covariates are included in the linear regression model or not. In a Bayesian framework, and for a judicious choice of prior distributions, we can explicitly calculate the posterior distribution $\pi$ up to a constant.

We want to sample from this distribution in order to approximate quantities like the expected value $\evplus{\pi}{\v\gamma}$, that is the marginal probability of inclusion of each variable. Often, the marginal probabilities provide a richer picture of the posterior distribution than a collection of modes found using stochastic optimisation techniques.

\subsection{Global versus local methods}
Our Sequential Monte Carlo approach to variable selection views a well studied problem from a different angle and provides new perspectives. The reason is two-fold.

Firstly, there is growing evidence that global methods, which track a population of particles, initially well spread over the sampling space $\B^d$, are often more robust than local methods based on Markov chain Monte Carlo. The latter are more prone to get trapped in the neighbourhood of local modes. We largely illustrate this effect in our simulations in Section \ref{sec:simu}.

Secondly, global methods have the property to be easily parallelisable. Parallel implementations of Monte Carlo algorithms have gained a tremendous interest in the very recent years \citep{lee2010utility,suchard2010some}, due to the increasing availability of multi-core (central or graphical) processing units in standard computers.

\subsection{Plan and notations}
The paper is organised as follows.

In Section \ref{sec:var sel}, we recapitulate the basics of Bayesian variable selection in linear regression models as the motivating application.

In Section \ref{sec:mcmc}, we briefly review the principal Markov chain Monte Carlo methods which are commonly used to integrate with respect to a binary distributions.

In Section \ref{sec:smc}, we describe an alternative approach to the same problem using Sequential Monte Carlo methods. The key ingredient of this algorithm is a parametric family which is flexible enough to come close to the target distribution.

In Section \ref{sec:bv}, we extensively discuss approaches for constructing rich parametric families on binary spaces. This is the core of our work. Some of the binary models discussed are not suitable in the framework of our Sequential Monte Carlo algorithm but mentioned for completeness of the survey.

In Section \ref{sec:simu}, we construct two examples of variable selection problems which yield challenging posterior distributions. We show that standard Markov chain techniques fail to produce reliable estimates of the marginal probabilities while our Sequential Monte Carlo approach successfully copes with the integration problem.

\paragraph{Notation}
For a vector $\v x\in \mathcal{X}^d$, we write $\v x_M$ for the sub-vector indexed by $M\subseteq\lbrace 1,\dots,d\rbrace$. We write $\v x_{i:j}$ if the indices are a complete sequence $i,\dots,j$. We denote by $\v x_{-i}$ the sub-vector $\v x_{\lbrace 1,\dots,d\rbrace\setminus\lbrace i\rbrace}$. We write $\abs{\v x}$ for $\sum_{k=1}^d \v x_k$.

For a matrix $\m A$, we denote its components by $a_{ij}$, its determinant by $\abs{\m A}$. The operator $\diag{\v x}$ transforms the vector $\v x$ into a diagonal matrix. For a finite set $M$, we denote by $\# M$ the number of elements in $M$.

%

%


%

%

\section{Variable selection: A binary sampling problem}
\label{sec:var sel}
The standard linear normal model postulates that the relationship between the observed explained variable $\v{y}\in\R^m$ and the observations $\m{Z}=[\v z_1,\dots,\v z_d]\in\R^{m, d}$ is
\begin{equation*}
\v y \mid \v\beta,\v\gamma,\sigma^2, \m Z \sim \normal{\m{Z}\,\diag{\v\gamma}\v{\beta}}{\sigma^2\m{I}_m}.
\end{equation*}
Here, $\v \beta$ is a vector of regression coefficients and $\sigma^2$ the variance of $\v y$. We denote by $\m{I}_m$ the identity matrix and assume the first column $\m Z_{\cdot,1}$ to be constant. The parameter $\v\gamma\in\B^d=\lbrace0,1\rbrace^d$ determines which covariates are included in or dropped from the linear regression model. In total, we can construct $2^d$ different linear normal models from the data.

We assign a prior distribution $\pi(\v\beta,\sigma^2,\v\gamma\mid\m Z)$ to the parameters. From the posterior distribution
\begin{equation*}
\pi(\v\beta,\sigma^2,\v\gamma\mid\v y,\m Z)\propto\pi(\v
y\mid\v\beta,\sigma^2,\v\gamma,\m Z)\,
\pi(\v\beta,\sigma^2,\v\gamma\mid\m Z) 
\end{equation*}
we may compute the posterior probability of each model
\begin{equation}
\label{eqn:int posterior}
\pi(\v\gamma\mid\v y,\m Z)=\int \pi(\v\beta,\sigma^2,\v\gamma\mid\v y,\m Z)\,d(\v\beta,\sigma^2)
\end{equation}
by integrating out the parameters $\v\beta$ and $\sigma^2$.

\paragraph{Hierarchical Bayesian model}
\label{sec:hb}
In a purely Bayesian context, we obtain, up to a constant, an explicit formula for the integral in \eqref{eqn:int posterior} by choosing conjugate hierarchical priors, that is a normal $\pi(\v{\beta}\mid\sigma^2,\v \gamma,\m Z)$ and an inverse-gamma $\pi(\sigma^2\mid\v\gamma,\m Z)$. For all Bayesian posterior distributions in this paper, we use the prior distributions
\begin{align*}
&\pi(\v{\beta}\mid\sigma,\v \gamma,\m Z)
=\normal{\v{0}}{\sigma^2 v^2 \diag{\v \gamma}}, & \sigma^2>0, \\
&\pi(\sigma^2\mid\v\gamma,\m Z)
=\mathcal{I}(w/2,\lambda w/2), & w >0,\ \lambda>0, \\
&\pi(\v\gamma\mid\m Z)=\mathcal{U}(\B^d),
\end{align*}
where $\mathcal{I}$ denote an Inverse-Gamma and $\mathcal{U}$ a uniform law.

For our numerical examples in Section \ref{sec:simu}, we assume not to have any prior information about the data. We follow the recommendations of \citet{george_mcculloch_97} and choose the hyper-parameters
\begin{equation}
w=4.0,\quad\lambda=\widehat\sigma^2_{\v 1},\quad v^2=10.0/\lambda,
\end{equation}
where $\widehat\sigma^2_{\v 1}$ is the least square estimate of $\sigma^2$ based on the saturated model. The rationale behind this choice is to have a flat prior on $\v \beta$ and provide $\sigma^2$ with sufficient mass on the interval $(\hat\sigma^2_{\v 1},\hat\sigma^2_{\v 0})$, where $\hat\sigma^2_{\v 0}$ denotes the variance of $\v y$.

Next, we quickly state the form of the log-posterior mass function. We write $\m{Z}_{\v\gamma}$ for $\m{Z}\,\diag{\v\gamma}$ without zero columns. Let $\v b_{\v\gamma}=\m{Z}_{\v\gamma}\t\,\v{y}$ and
\begin{equation}
\label{eq:chol}
\m C_{\v\gamma,v}^{} \m C_{\v\gamma,v}\t = \m{Z}_{\v\gamma}\t\,\m{Z}_{\v\gamma}^{}+v^{-2}\m I_{\abs{\v\gamma}} 
\end{equation}
a Cholesky decomposition. We denote the least square estimate of $\sigma^2$ based on $\nu$ and the model $\v\gamma$ by
\begin{equation*}
\widehat\sigma^2_{\v\gamma,v}=\frac{1}{m}\left(\v{y}\t\v{y}-(\m C^{-1}_{\v\gamma,v}\v b_{\v\gamma}^{})\t(\m C^{-1}_{\v\gamma,v}\v b_{\v\gamma}^{})\right).
\end{equation*}
We find the log-posterior probability to be
\begin{align*}
\log\pi(\v\gamma\mid\v y,\m Z)
&=\, \mu -\textstyle\sum_{i=1}^{\abs{\v \gamma}} \log c_{i,i}^{(\v\gamma,v)}-\abs{\v\gamma}\log(v) \\ &\quad-\frac{w+m}{2}\log(w\lambda/m+\widehat\sigma^2_{\v\gamma,v}),
\end{align*}
where $\mu$ is an unknown normalization constant.

\paragraph{Related approaches}
In a Frequentist framework, we choose a model which minimizes some specified criterion. A popular one is Schwarz's Criterion \citep[also Bayesian Information Criterion]{schwarz_78} which basically is a second degree Laplace approximation of
\eqref{eqn:int posterior}:
\begin{equation*}
\log\pi(\v\gamma\mid\v y,\m Z)\approx \mu
-\frac{\abs{\v\gamma}}{2}\log(m) -\frac{m}{2}\log(\widehat\sigma^2_{\v\gamma}),
\end{equation*}
where $\widehat\sigma_{\v \gamma}^2=\lim_{v\to\infty}\widehat\sigma_{\v \gamma,v}^2$ is the maximum likelihood estimator of $\sigma^2$ based on the model $\v \gamma$. Note that for a large sample size $m$ the Hierarchical Bayesian approach and the Bayesian Information Criterion coincide.

\paragraph{Alternative approaches}
The posterior of a Bayesian linear regression variable selection problem has, in general, no particular structure we can exploit to speed up optimisation or integration with respect to $\pi$. Therefore, alternative approaches such as the Least Absolute Shrinkage and Selection Operator \citep{tibshirani1996regression} have been proposed which draw from the theory of convex optimization and allow for computation of larger problems.

While a comparison between competing approaches to variable selection is beyond the scope of this paper, we remark that more sophisticated, parallelisable algorithms are essential for making Bayesian modelling feasible in the context of high dimensional problems where alternative methods are often used for practical reasons only.

%

%


%

%

\section{Markov chain Monte Carlo on binary spaces}
\label{sec:mcmc}
Markov chain Monte Carlo is a well-studied approach to approximate the expected value of a posterior $\pi$ given by a Bayesian model choice problem \citep{george_mcculloch_97}. In this section, we rapidly review the standard methods we are going to compare our Sequential Monte Carlo approach against.

There are more advanced Markov chain Monte Carlo algorithms that use parallel tempering ideas combined with more elaborate local moves \citep[see e.g.][]{liang2000evolutionary, bottolo2010ess}, but a thorough comparison is beyond the scope of this paper. For background on Markov chain Monte Carlo methods, we refer to standard literature \citep[see e.g.][chaps. 7-12]{RobCas}.

\subsection{Framework}
The idea is to construct a transition kernel $\kappa$, typically some version of a Metropolis-Hastings kernel, which admits $\pi$ as unique stationary distribution. Then, the distribution of the Markov chain $\v x_{t+1}\sim\kappa(\v x_t,\cdot)$ started at some randomly chosen point $\v x_0\in\B^d$ converges to $\pi$.

We obtain an estimate $\evplus{\pi}{\v\gamma}\approx n^{-1}\sum_{t=b}^{n+b}\v x_t$ for the expected value via the ergodic theorems for Markov chains. The first $b$ states are usually discarded to give the chain some time to converge towards the stationary distribution before we start to average. For the estimate to be valid, we need to ensure that the at time $b$ the chain is close to its stationary distribution $\pi$, and at time $n+b$ we have sampled an ergodic trajectory such that the ergodic theorems applies.

Classic Markov chain methods on binary spaces work locally, that is they propose moves to neighbouring models in the Metropolis-Hastings steps. A neighbouring model is a copy of the current model where just a few components are altered. We shall see that these kinds of transition kernels often fail to sample ergodic trajectories within a reasonable amount of time if the stationary distribution $\pi$ is very multi-modal.

\paragraph{Algorithm}
We loop over a uniformly drawn subset of components $I\sim\mathcal{U}(\lbrace M\subseteq \lbrace1,\dots,d\rbrace\mid \#M=k\rbrace)$ and propose to change the components $i\in I$. The number of components $k$ might be fixed or drawn from some distribution $\mathcal{G}$ on the index set $\lbrace1,\dots,d\rbrace$.

Precisely, we take a copy $\v y$ of the current state $\v x_t$ and replace $y_i$ by $Y_i\sim\ber_{p_i(\v x)}$ for all $i\in I$, where
\begin{equation*}
\ber_{p_i(\v x)}(\gamma)=p_i(\v x)^{\gamma_i}(1-p_i(\v x))^{1-\gamma_i}
\end{equation*}
is a Bernoulli distribution with parameter $p_i(\v x)\in(0,1)$. We set $\v x_{t+1}=\v y$ with probability
\begin{equation}
\label{eq:acc ratio}
\frac{\pi(\v y)}{\pi(\v x_t)}
\frac{\prod_{i\in I} \ber_{p_i(\v y)}(\v x_t)}{\prod_{i\in I} \ber_{p_i(\v x_t)}(\v y)}\wedge 1,
\end{equation}
and $\v x_{t+1}=\v x_t$ otherwise. This framework, summarized in Algorithm \ref{algo:generic mg}, yields a Markov chain with unique invariant distribution $\pi$ for any fixed $\v p\in(0,1)^d$. The interesting special cases, however, use a $p(x)$ which depends on the current state of the chain.

\floatname{algorithm}{Algorithm}
\begin{algorithm}[H]
\caption{Generic metropolised Gibbs kernel}
\label{algo:generic mg}
\begin{algorithmic}
\onehalfspacing
\REQUIRE $\v x\in\B^d$
\STATE $U\sim\mathcal{U}([0,1]),\ k\sim \mathcal{G}_{k^*}$
\STATE $I\sim\mathcal{U}(\lbrace M\subseteq \lbrace1,\dots,d\rbrace\mid \#M=k\rbrace)$
\STATE $\v y\gets \v x$
\STATE \textbf{for} $i\in I$ \textbf{ do } $y_i \sim \ber_{p_i(\v x)}$ \\[0.5em]
\STATE \textbf{if} {$\displaystyle \,\frac{\pi(\v y)}{\pi(\v x)}
\frac{\prod_{i\in I} \ber_{p_i(\v y)}(\v x)}{\prod_{i\in I} \ber_{p_i(\v x)}(\v y)}
> U$} \textbf{ then } $\v x \gets \v y$\\[0.5em]
\RETURN $\v x$
\end{algorithmic}
\end{algorithm}

\paragraph{Performance}
We refer to the ratio \eqref{eq:acc ratio} as the acceptance pro\-bability of the Metropolis-Hastings step. In binary spaces, however, accepting a proposal does not imply we are changing the state of the chain, since we are likely to re-propose the current state $y=\v x_t$. We are actually interested in how fast the chain explores the state spaces, precisely its mutation probability $\prob{\v x_{t+1} \neq \v x_t}$.

\subsection{Standard Markov chain methods}
For this section, let $k=1$ be constant. Algorithm \ref{algo:generic mg} collapses to changing a single component. Instead of independently drawing the index $i\sim\mathcal{U}(\lbrace1,\dots,d\rbrace)$, we could also iterate $i$ through a uniformly drawn permutations $\sigma(\lbrace1,\dots,d\rbrace)$ of the index set $\lbrace1,\dots,d\rbrace$.

Kernels of this kind are often referred to as metropolised Gibbs samplers, since they proceed component-wise as does the classical Gibbs sampler, but also involve a Metropolis\-Hastings step. In the sequel, we discuss some special cases.

\paragraph{Classical Gibbs} The Gibbs sampler sequentially draws each component from the full marginal distribution, which corresponds to
\begin{align*}
p_i(\v x)
&\eqdef\pi(\gamma_i=1\mid \v\gamma_{-i}=\v x_{-i}) \\
&=\frac{\pi(\gamma_i=1,\v\gamma_{-i}=\v x_{-i})}{\pi(\gamma_i=1,\v\gamma_{-i}=\v x_{-i})+\pi(\gamma_i=0,\v\gamma_{-i}=\v x_{-i})}.
\end{align*}
By construction, the acceptance probability is $1$ while the mutation probability is only $\pi(\v y)/(\pi(\v x_t)+\pi(\v y))$, where $\v y$ is a copy of the current state $\v x_t$ with component $i$ altered.

\paragraph{Adaptive metropolised Gibbs}
An adaptive extension of the metropolised Gibbs has been proposed by \citet{nott2005adaptive}. The full marginal distribution $\pi(\gamma_j=1\mid \v\gamma_{-j}=x_{-j})$ is approximated by a linear predictor. In their notation,
\begin{align*}
p_i(\v x)\eqdef\left\lbrack \left(\psi_i-\frac{\m W_{-i}\v x_{-i}}{w_{i,i}}\right) \vee \delta \right\rbrack \wedge (1-\delta),
\end{align*}
where $\psi$ is the estimated mean, $\m W^{-1}$ the estimated covariance matrix and $\delta\in(0,1/2)$ a design parameter which ensures that $p_i(\v x)$ is a probability. Analogously to our vector notation, $\m W_{-i}$ denotes the matrix $\m W$ without the $i$th row and column. We obtain the estimates from the past trajectory of the chain $\v x_b,\dots,\v x_{t-1}$ and update them periodically.

The mutation probability is of the same order as that of the Gibbs kernel, but adaption largely avoids the computationally expensive evaluations of $\pi$. The non-adaptive Gibbs sampler already requires evaluation of $\pi(\v y)$ just to produce the proposal $\v y$. In contrast, the adaptive metropolised Gibbs samples from a proxy and only evaluates $\pi(\v y)$ if $\v y\neq \v x_t$.

\paragraph{Modified metropolised Gibbs}
In comparison to the classical Gibbs kernel, we obtain a more efficient chain \citep{liu1996peskun} using the simple form
\begin{equation*}
p_i(\v x)\eqdef1-x_i.
\end{equation*}
Since we always propose to change the current state, the acceptance and mutation probabilities are the same. They amount to $\pi(\v y)/\pi(\v x)\wedge1$, where $\v y$ is a copy of the current state $\v x$ with component $i$ altered. Comparing the mutation probabilities of the two kernels, we see that the modified metropolised Gibbs chain moves, on average, faster than the classical Gibbs chain.

\subsection{Block updating}
\label{sec:blockup}
The modified metropolised Gibbs easily generalises to the case where $k$ may take values larger than one. Suppose, for example, we propose to change
\begin{equation*}
k\sim \mathcal{G}_{k^*}(k) \propto \frac{(1-1/k^*)^{k-1}}{k^*}\,\ind_{\lbrace1,\dots,d\rbrace}(k)
\end{equation*}
components simultaneously, where $\mathcal{G}_{k^*}$ is a truncated geometric distribution. Note that we suggest, on average, to change approximately $k^*$ components. In other words, for larger values of $k^*$, we are more likely to propose further steps in the sampling space.

Large step proposals improve the mixing properties of the chain and help to escape from the attraction of local modes. They are, however, less likely to be accepted than single component steps which leads to a problem-dependent trade-off. In our numerical examples, we could not observe any benefit from block updating, and we do not further consider it to keep the comparison with our Sequential Monte Carlo method more concise.

\subsection{Independent proposals}
\label{sec:ind mh}
We can construct a fast mixing Markov chain based on independent proposals. Let $q$ be some distribution with $\pi \ll q$, that is $q(\v\gamma)=0\ \Rightarrow\ \pi(\v\gamma)=0$ for all $\v\gamma\in\B^d$. We propose a new state $y\sim q$ and accept it with probability
\begin{equation}
\label{eq:ind mh}
\frac{\pi(\v y)}{\pi(\v x_t)}\frac{q(\v x_t)}{q(\v y)}\wedge 1.
\end{equation}
The associated Markov chain has the unique invariant measure $\pi$. This kernel is referred to as the independent Metro\-polis-Hastings kernel, since the proposal distribution is not a function of the current state $\v x_t$. The mutation rate is the acceptance rate minus $q(\v x_t)$, so the two notions practically coincide in large spaces.

Obviously, in order to make this approach work, we need to choose $q$ sufficiently close to $\pi$, which implies high acceptance rates on average. In absence of reliable prior information, however, we are not able to produce such a distribution $q$. We shall, however, use precisely this Markov kernel as part of our Sequential Monte Carlo algorithm. In this context, we can calibrate sequences $q_t$ of proposal distributions to be close to our current particle approximation.

%

%


%

%

\section{Sequential Monte Carlo on binary spaces}
\label{sec:smc}
In this section, we show how to estimate the expected value with respect to a probability mass function $\pi(\v \gamma)$ defined on $\B^d$ using Sequential Monte Carlo \citep{del2006sequential}. This general class of algorithms alternates importance sampling steps, resampling steps and Markov chain transitions, to recursively approximate a sequence of distributions, using a set of weighted `particles' which represent the current distribution. In the following, we present a version which is tailored to work on binary spaces.

For readers not familiar with Sequential Monte Carlo, the following algorithm described might seem rather complex at first glance. We introduce the steps separately before we look at the complete algorithm. We give comprehensive instructions which correspond exactly to our implementation in order to make our results plausible and easily reproducible for the reader. 

\floatname{algorithm}{Procedure}
\setcounter{algorithm}{0}

\subsection{Building a sequence of distributions}
The first ingredient of Sequential Monte Carlo is a smooth sequence of distributions $(\pi_t)_{t=0}^\tau$, which ends up at the distribution of interest $\pi_\tau=\pi$. The intermediary distributions $\pi_t$ are purely instrumental: the idea is to depart from a distribution $\pi_0$ with broad support and to progress smoothly towards the distribution of interest $\pi$.

\paragraph{Initial distribution}
Theoretically, we can use any $\pi_0$ with $\pi \ll \pi_0$ that can sample from as initial distribution. Numerical experiments taught us, however, that premature adjustment of $\pi_0$, for example using Markov chain pilot runs, leads to faster but less robust algorithms.

Thus, in practice, we recommend the uniform distribution for its simplicity and reliability. Therefore, in the sequel, we let $\pi_0=\mathcal U(\B^d)$.

\paragraph{Intermediate distributions}
We construct a smooth sequence of distributions by judicious choice of an associated real sequence $(\varrho_t)_{t=0}^\tau$ increasing from zero to one. The most convenient and somewhat natural strategy is the geometric bridge \citep{gelman1998simulating,Neal:AIS,del2006sequential}
\begin{equation}
\label{eqn:geo bridge}
\pi_t(\v\gamma)\ \proptodef\ \pi_0(\v\gamma)^{1-\varrho_t}\pi(\v\gamma)^{\varrho_t}\propto \pi(\v\gamma)^{\varrho_t}.
\end{equation}
Alternatively, one could use a sequences of mixtures
\begin{equation*}
\pi_t(\v\gamma)\proptodef (1-\varrho_t)\pi_0(\v\gamma)+\varrho_t\pi(\v\gamma)
\end{equation*}
or, in a Bayesian context, a sequences of posterior distributions where data is added as $\varrho_t$ increases, that is
\begin{equation*}
\pi_t(\v\gamma)=\pi(\v\gamma\mid \v z_1,\dots,\v z_{\lfloor\varrho_t m\rfloor}),
\end{equation*}
see \citep{chopin2002sequential}. In the following, we use the geometric bridge \eqref{eqn:geo bridge} for its computational simplicity and present a procedure to determine a suitable sequence $(\varrho_t)_{t=0}^\tau$.

\subsection{Assigning importance weights}
Suppose we have already produced a sample $\v x_1\up{\,t-1},\dots,\v x_n\up{\,t-1}$ of size $n$ from $\pi_{t-1}$. We can roughly approximate $\pi_t$ by the empirical distribution
\begin{equation}
\label{eq:emp approx}
\pi_t(\v\gamma)\approx\sum_{k=1}^n w_t(\v x_k^{[\,t-1]})\,\delta_{\v x_k^{[\,t-1]}}(\v\gamma),
\end{equation}
where the corresponding importance function $w_t$ is
\begin{equation}
\label{eq:imp weights}
w_t(\v x_k)\eqdef\frac{u_t(\v x_k)}{\sum_{i=1}^n u_t(\v x_i)},\quad u_t(\v x)\eqdef\frac{\pi_t(\v x)}{\pi_{t-1}(\v x)}=\pi^{\alpha_t}(\v x).
\end{equation}

Note that $\alpha_t=\varrho_{t}-\varrho_{t-1}$ is the step length at time $t$. As we choose $\alpha_t$ larger, that is $\pi_t$ further from $\pi_{t-1}$, the weights become more uneven and the accuracy of the importance approximation deteriorates.

\begin{algorithm}
\caption{Importance weights}
\label{algo:weight}
\begin{algorithmic}
\REQUIRE{$\alpha,\ \pi,\ \m X=(\v x_1,\dots,\v x_n)\t$}
\STATE $u_k\gets\pi^{\alpha}(\v x_k)$ \textbf{ for all } $k=1,\dots,n$
\STATE $w_k\gets u_k/(\sum_{i=1}^nu_i)$ \textbf{ for all } $k=1,\dots,n$
\RETURN $\v w=(w_1,\dots,w_n)$
\end{algorithmic}
\end{algorithm}

If we repeat the weighting steps until we reach $\pi_\tau=\pi$, we obtain a classical importance sampling estimate with instrumental distribution $\pi_0$. The idea of the  Sequential Monte Carlo algorithm, however, is to control the weight degeneracy such that we can intersperse resample and move steps before loosing track of our particle approximation.

\paragraph{Effective sample size}
We measure the weight degeneracy through the effective sample size criterion, see \citep{KongLiuWong}. In our case, we have
\begin{equation*}
\eta(\alpha,\m X) \eqdef \frac{\left(\sum_{k=1}^n w_\alpha(\v x_k)\right)^2}{n\sum_{k=1}^n w_\alpha(\v x_k)^2}
=\frac{\left(\sum_{k=1}^n \pi^{\alpha}(\v x_k)\right)^2}{n\sum_{k=1}^n \pi^{\alpha}(\v x_k)^2}\in[1/n,1].
\end{equation*}
The effective sample size is $1$ if all weights are equal and $1/n$ if all mass is concentrated in a single particle.

For a geometric bridge \eqref{eqn:geo bridge}, the effective sample size is merely a function of $\alpha$. We can thus control the weight degeneracy by judicious choice of the step lengths $\alpha_t$.

\subsection{Finding the step length}
We pick a step length $\alpha$ such that the effective sample size $\eta(\alpha)$ equals a fixed value $\eta^*$, see \citep{jasra2008inference}. Since $\eta$ is continuous and monotonously increasing in $\alpha$, we solve 
\begin{equation}
\label{eq:ess}
\eta(\alpha,\m X)=\eta^*
\end{equation}
using bi-sectional search, see Procedure \ref{algo:step}. This approach is numerically more stable than a Newton-Raphson iteration, for the derivative $\partial \eta(\alpha,\v x)/\partial \alpha$ involves fractions of sums of exponentials which are difficult to handle.

Let $\alpha^*$ be the unique solution to \eqref{eq:ess}. We can construct an associated sequence setting $\varrho_t=1 \wedge (\varrho_{t-1}+ \alpha^*)$. Thus, the number of steps $\tau$ depends on the complexity of the integration problem at hand and is not known in advance.

In other words, for fixed $\eta^*$, the associated sequence $(\varrho_t)_t^\tau$ is a self-tuning parameter. In our simulations, we always choose $\eta^*=0.9$, which yields convincing results on both example problems in Section \ref{sec:simu}. Smaller values significantly speed up the Sequential Monte Carlo algorithm but lead to a higher variation in the results. 

\begin{algorithm}
\caption{Find step length}
\label{algo:step}
\begin{algorithmic}
\REQUIRE{$\varrho,\,\m X=(\v x_1,\dots,\v x_n)\t$}
\STATE $l\gets 0,\, u \gets 1.05-\rho,\,\alpha\gets0.05$
\REPEAT
	\STATE \textbf{if} {$\eta(\alpha,\m X)  < \eta^*$} \textbf{then} $u \gets \alpha,\, \alpha \gets (\alpha + l)/2$
	\STATE \textbf{else} $l \gets \alpha,\, \alpha \gets (\alpha + u)/2$
\UNTIL{$\abs{u-l}<\varepsilon$ \OR $l > 1-\varrho$}
\RETURN $\alpha\wedge(1-\varrho)$
\end{algorithmic}
\end{algorithm}

\subsection{Resampling the system}
Suppose we have a sample $\m X\up{t-1}=(\v x_1\up{t-1},\dots,\v x_n\up{t-1})$ of size $n$ from $\pi_{t-1}$ with importance weights as defined in \eqref{eq:imp weights}. We can obtain a sample $\widehat{\m X}\up{t}=(\hat{\v x}_1\up{t},\dots,\hat{\v x}_n\up{t})$ which is approximately distributed according to $\pi_t$ by drawing from the empirical approximation defined in \eqref{eq:emp approx}.

\begin{algorithm}
\caption{Resample (systematic)}
\label{algo:resample}
\begin{algorithmic}
\REQUIRE{$\v w=(w_1,\dots,w_n),\,\m X=(\v x_1,\dots,\v x_n)\t$}
\STATE $v\gets n\,w,\ j\gets1,\ c\gets v_1$
\STATE \textbf{sample} $u\sim\mathcal{U}([0,1])$
\FOR{$k=1,\dots,n$}
	\WHILE{$c < u$}
		\STATE $j\gets j+1,\ c \gets c+v_j$
        \ENDWHILE
        \STATE $\hat{\v x}_k\gets \v x_j,\ u\gets u+1$
\ENDFOR
\RETURN $\widehat{\m X}=(\hat{\v x}_1\dots,\hat{\v x}_n)\t$
\end{algorithmic}
\end{algorithm}

For the implementation of the resampling step, there exist several recipes. We could apply a multinomial resampling \citep{Gordon} which is straightforward. There are, however, more efficient ways like residual \citep{LiuChen}, stratified \citep{kitagawa1996monte} and systematic resampling \citep{CarClifFearn}. We use the latest in our simulations, see Procedure \ref{algo:resample}.

In the resulting unweighted particle system $\widehat{\m X}\up{t}$ of size $n$, the particles with small weights have vanished while the particles with large weights have bee multiplied. There are approaches that resample a weighted particle system of size $n$ from an augmented system of size $m>n$, see \citep{fearnhead2003line}, but these techniques are computationally more demanding without visibly improving our numerical results. Theoretically, one would expect a Rao-Blackwellisation effect but its analysis is beyond the scope of this paper.

In any case, if we repeat the weighting and resampling steps several times, we rapidly deplete our particle reservoir reducing the number of different particles to a very few. Thus, the particle approximation will be totally inaccurate. The key to fighting the decay of our approximation is the following move step.

\subsection{Moving the system}
\label{sec:move}
The resampling step provides an unweighted particle system $\widehat{\m X}\up{t}=(\hat{\v x}_1\up{t},\dots,\hat{\v x}_n\up{t})$ of $\pi_t$ containing multiple copies of many particles. The central idea of the Sequential Monte Carlo algorithm is to diversify the resampled system, replacing the particles by draws from a Markov kernel $\kappa_t$ with invariant measure $\pi_t$ \citep{GilksBerzu}.

Since the particle $\v x_k\up{0}$ is, approximately, distributed according to $\pi_t$, a draw $\v x_k\up{1}\sim\kappa_t(\v x_k\up{0},\cdot)$ is again, approximately, distributed according to $\pi_t$. We can repeat this procedure over and over without changing the target of the particle approximation.

Note that, even if the particles $\v x_{k}\up{0}=\dots=\v x_{m}\up{0}$ are equal after resampling, the particles $\v x_{k}\up{s},\dots,\v x_{m}\up{s}$ are almost independent after sufficiently many move steps. In order to make the algorithm practical, however, we need a transition kernel which is rapidly mixing and therefore diversifies the particle system within a few steps. Therefore, the locally operating Markov kernels reviewed in Section \ref{sec:mcmc} are not suitable. In fact, our numerical experiments suggest that making a Sequential Monte Carlo algorithm work with local kernels is practically impossible.

Therefore, we use a Metropolis-Hastings kernel with independent proposals as described in Section \ref{sec:ind mh}. Precisely, we construct a kernel $\kappa_t$ employing a parametric family $q_\theta$ on $\B^d$ which, for some $\theta$, is sufficiently close to $\pi_t$ to allow for high acceptance probabilities.

For this purpose, we fit a parameter $\theta_t$ to the particle approximation $(\v w_t,\m X_t)$ of $\pi_t$ according to some convenient criterion. The choice of the parametric family $q_\theta$ is crucial to a successful implementation of the Sequential Monte Carlo algorithm. We come back to this issue in Section \ref{sec:bv}.

\begin{algorithm}
\caption{Move}
\label{algo:move}
\begin{algorithmic}
\REQUIRE{\parbox{0.4\textwidth}{$\m X\up{0}=(\v x_1\up{0},\dots,\v x_n\up{0})$ \textbf{ targeting } $\pi_t$ \\[0.2em]
$\kappa(\v y,\gamma)$ such that $\pi_t(\gamma)=\sum_{y\in\B^d} \pi_t(\v y) \kappa(\v y,\gamma)$}}
\STATE $s\gets 1$
\REPEAT
	\STATE \textbf{sample} ${\v x}_k\up{s}\sim \kappa ({\v x}_k\up{s-1},\cdot)$ \textbf{for all} $k=1,\dots,n$ \\[0.5em]
\UNTIL{$|\zeta(\m X\up{s})-\zeta(\m X\up{s-1})|<0.02$ \OR $\zeta(\m X\up{s})>0.95$} \\[0.5em]
\RETURN $\m X\up{s}=(\v x_1\up{s}\dots,\v x_n\up{s})\t$
\end{algorithmic}
\end{algorithm}

\paragraph{Particle diversity}
We need to determine how often we move the particle system before we return to the weight-resample step. An easy criterion for the health of the particle approximation $\m X=(x_1,\dots,x_n)$ is its particle diversity
\begin{equation}
\label{eqn:pd}
\zeta(\m X) \eqdef \frac{\#\lbrace \v x_k\mid k=1,\dots,n\rbrace}{n}\in[1/n,1],
\end{equation}
that is the proportion of distinct particles. Note that the particle diversity is a quality criterion which has no simple analogue in continuous sampling spaces.

For optimal results, we recommend to keep on moving the particle system until the particle diversity cannot be augmented any longer. In the first steps of the algorithm, $\pi_t$ is still close to the uniform distribution, and we manage to raise the particle diversity up to one. As $\pi_t$ is approaching a strongly multi-modal target distribution $\pi$, however, the particle diversity reaches a steady-state we cannot push it beyond.

Clearly, even if we could draw a particle system independently from $\pi$, the particle diversity would be a lot smaller than one, since we would draw the modes of $\pi$ several times.

\paragraph{Aggregated weights}
Shifting weights between identical particles does not affect the nature of the approximation but it changes the effective sample size $\eta(\v w)$ which seems paradoxical at first sight. For reasons of parsimoniousness, we could just keep a single representative $\v x_*$ for identical particles $\v x_{*_1}=\dots=\v x_{*_k}$ and aggregate the associated weights to the sum $w_*=w_{*_1}+\cdots+w_{*_k}$ without changing the quality of the particle approximation. There are, however, three reasons why we refrain from doing so.

Firstly, it is vital not to confuse the weight disparity induced by reweighting according to the progression of $\pi_t$ and the weight disparity due to aggregation of the weights of multiply sampled states. From the aggregated system, we cannot tell whether the effective sample size is determined by the gap between $\pi_t$ and $\pi_{t+1}$, that is the step length $\alpha$, or by the presence of particle copies due to the mass of $\pi_t$ being very concentrated. Therefore, it seems more difficult to control the smoothness of the sequence of distributions and find a suitable sequence $(\varrho_t)_{t=0}^\tau$.

Secondly, aggregation is an additional computational effort equivalent to keeping the particle system sorted. Here, we trade in computing time for memory while the required memory is proportional to the number of particles and not critical in the context of our algorithm.

Thirdly, the straightforward way to implement repeated move steps is breaking up the particles into multiple copies corresponding to their weights and moving them separately. Consequently, instead of permanently splitting and aggregating the weights we might just allow for multiple copies of the particles.

\subsection{The Resample-move algorithm}
Finally, we summarize the complete Sequential Monte Carlo method in Algorithm \ref{algo:smc}. Note that, in practice, the sequence $\pi_t=\pi^{\rho_t}$ is not indexed by $t$ but rather by $\rho_t$, that is the counter $t$ is only given implicitly.

\floatname{algorithm}{Algorithm}
\setcounter{algorithm}{1}
\begin{algorithm}
\caption{Resample-move}
\label{algo:smc}
\begin{algorithmic}
\REQUIRE{$\pi\colon \B^d\to[0,\infty)$}
\STATE \textbf{sample} $\v x_k\stackrel{\mathrm{iid}}{\sim}\mathcal{U}(\B^d)$ \textbf{for all} $k=1,\dots,n$.
\STATE \vspace{0.5em}
\begin{tabular}{lll}
$\alpha	$&\hspace{-3mm}$\gets\textbf{find step length}(0, \m X)$ & (Procedure \ref{algo:step}) \\[0.5em]
$\v w	$&\hspace{-3mm}$\gets\textbf{importance weights}(\alpha,\pi,\m X)$ & (Procedure \ref{algo:weight}) \\[0.5em]
\end{tabular}
\WHILE {$\varrho<1$}
\STATE \vspace{0.5em}
\begin{tabular}{lll}
$q_\theta	$&\hspace{-3mm}$\gets\textbf{fit binary model}(w, \m X)$ 		&\hspace{-3mm} (Section \ref{sec:bv}) \\[0.5em]
$\widehat{\m X}	$&\hspace{-3mm}$\gets\textbf{resample}(w,\m X)$ 			&\hspace{-3mm} (Procedure \ref{algo:resample}) \\[0.5em]
$\m X		$&\hspace{-3mm}$\gets\textbf{move}(\kappa_{\pi, q_\theta}, \widehat{\m X})$&\hspace{-3mm} (Procedure \ref{algo:move}) \\[0.5em]
$\alpha		$&\hspace{-3mm}$\gets\textbf{find step length}(\rho, \m X)$		&\hspace{-3mm} (Procedure \ref{algo:step}) \\[0.5em]
$\v w		$&\hspace{-3mm}$\gets\textbf{importance weights}(\alpha,\pi,\m X)$	&\hspace{-3mm} (Procedure \ref{algo:weight}) \\[0.5em]
$\rho		$&\hspace{-3mm}$\gets\varrho+\alpha$& \\[0.5em]
\end{tabular}
\ENDWHILE
\vspace{0.3em}
\RETURN $\sum_{k=1}^n w_k\v x_k\approx\evplus{\pi}{\v \gamma}$
\end{algorithmic}
\end{algorithm}

For an efficient implementation, we recommend to store the values $\pi(\v x_1),\dots,\pi(\v x_n)$ and $q_\theta(\v x_1),\dots,q_\theta(\v x_n)$ to avoid unnecessary evaluations. When updating the latter set, we can exploit the fact that, in a systematically resampled particle system, multiple copies of the same particles are neighbours.

%

%


%

%

\section{Multivariate binary models}
\label{sec:bv}
In this section, we address the choice of a multivariate binary parametric family $q_\theta$ with parameter $\theta\in\Theta$ needed to construct the independent Metropolis-Hastings kernel used in Procedure \ref{algo:move}. 

\subsection{Desired properties}
\label{sec:model prop}
We first frame the properties making a parametric family suitable for our Sequential Monte Carlo algorithm.
\begin{enumerate}[(a)]
\item For reasons of parsimony, we want to construct a family of distributions with at most $\mathrm{dim}(\theta)\leq d(d+1)/2$ parameters. More complex families are usually computationally too expensive to handle.
\item Given a sample $\m X=(\v x_1,\dots,\v x_n)$ from the target distribution $\pi$, we want to estimate $\theta^*$ such that the binary model $q_{\theta^*}$ is close to $\pi$. For instance, $\theta^*$ might be a maximum likelihood or method of moments estimator.
\item We want to generate independent samples from $q_\theta$. If we can compute the conditional or marginal distributions, we can write $q_\theta$ as
\begin{align}
\textstyle
\label{eqn:marg dec}
q_{\theta}(\v\gamma)
&=q_{\theta}(\v\gamma_1)\prod_{i=2}^d q_{\theta}(\v\gamma_i\vert\v\gamma_{1:i-1}) \\
&=q_{\theta}(\v\gamma_1)\prod_{i=2}^d q_{\theta}(\v\gamma_{1:i})/q_{\theta}(\v\gamma_{1:i-1}). \nonumber
\end{align}
Using the chain rule decomposition \eqref{eqn:marg dec}, we can sample a random vector $\v\gamma\sim q_{\v\theta}$ component-wise, conditioning on the entries we already generated.
\item We need to rapidly evaluate $q_\theta(\v\gamma)$ for any $\v\gamma\in\B^d$ in order to compute the Metropolis-Hastings ratio \eqref{eq:ind mh}.
\item Analogously to the multivariate normal, we want our calibrated binary model $q_{\theta^*}$ to produce samples with the mean and covariance of $\pi$. If $q_\theta$ is not flexible enough to capture the dependence structure of $\pi$, the Metropolis-Hastings kernel in Procedure \ref{algo:move} cannot provide satisfactory acceptance rates for complex target distributions $\pi$.
\end{enumerate}
In the following we construct a suitable parametric family and explain how to deploy it in Algorithm \ref{algo:smc}.

Most of the literature on binary data stems from response models, multi-way contingency tables and multivariate interaction theory \citep{cox1972analysis}. For completeness, we append a brief list of other binary models mentioned in the literature which fail, for various reasons, to work in Sequential Monte Carlo applications. Providing parametric families which meet the above requirements in high dimensions is a difficult task and understanding the shortcomings of alternative approaches an important part of the discussion.

\subsection{The logistic conditionals model}
\label{sec:logistic model}
In the previous paragraph, we already mentioned that a factorization \eqref{eqn:marg dec} of the mass function $q_\theta(\v\gamma)$ into conditional distributions permits to sample from the parametric family. Unfortunately, for a complex $d$-dimensional binary model, we usually cannot calculate closed-form expressions for the conditional or marginal mass functions.

We get around the computation of the marginal distributions of $q_\theta(\v\gamma)$ if we directly fit univariate models $q_{\v b_i}(\gamma_i\mid\gamma_{1:i-1})$ to the conditionals $\pi(\gamma_i\mid\gamma_{1:i-1})$ of the target function. \citet{qaqish2003family} suggested the use of linear regressions to model the conditional probabilities. This approach, however, does not guarantee that the fitted model is a valid distribution since the mass function might be negative.

\paragraph{Construction of the model}
We propose to rather use logistic regressions for the conditional probabilities. Precisely, we adjust the univariate models
\begin{equation*}
\textstyle
\logit(\probplus{\pi}{\gamma_i=1\mid\v\gamma_{1:i-1}})\eqdef b_{i,i}+\sum_{j=1}^{i-1} b_{i,j} \gamma_j,\quad i=1,\dots,d
\end{equation*}
where $\logit(p)=\log p- \log(1-p)$. In the context of our Sequential Monte Carlo application, we take the particle system $\m X$ and regress $\v y\up{i}=\m X_i$ on the columns $\m Z\up i=(\m X_{1:i-1}, \v 1)$, where the column $\v Z_i\up i$ yields the intercept to complete the logistic model.

For a $d$-dimensional lower triangular matrix $\m B$, we define the logistic conditionals model as
\begin{align}
\label{eqn:lb}
\textstyle
q_{\m B}(\v \gamma)
\eqdef \prod_{i=1}^d
\ber_{p(b_{i,i}+\v b_{i,1:i-1} \v \gamma_{1:i-1}\t)}(\gamma_i)
\end{align}
where $p(y)=\logit^{-1}(y)=(1+\exp(-y))^{-1}$ is the logistic function. Recall that $\ber_p(\v\gamma)=p^{\gamma}(1-p)^{1-\gamma}$ denotes the univariate Bernoulli distribution with expected value $p\in[0,1]$.

There are $d!$ possible logistic regressions models and we arbitrarily pick one while there should be a parametrization which is optimal in a sense of nearness to the data $\m Z$. We observed, however, that permuting the components had, in practice, no impact on the quality of the approximation.

Keep in mind that the number of observations in the logistic regressions is the size $n$ of the particle system and typically very large. For instance, we run our numerical examples in Section \ref{sec:simu} using $n=2\times10^4$ particles. Therefore, the fit of the logistic regressions is usually very good.

\paragraph{Sparse version}
The major drawback of any kind of multiplicative model is the fact that we have no closed-form likelihood-maximizers, and therefore the parameter estimation requires costly iterative fitting procedures. Therefore, even before discussing the fitting procedure, we construct a sparse version of the logistic conditionals model which we can estimate faster than the saturated model.

Instead of fitting the saturated model $q(\gamma_i\mid\gamma_{1:i-1})$, we preferably work with a more parsimonious regression model like $q(\gamma_i\mid\gamma_{L_i})$ for some index set $L_i\subseteq \lbrace1,\dots,i-1\rbrace$, where the number of predictors $\#L_i$ is typically smaller than $i-1$. We solve this nested variable selection problem using some simple, fast to compute criterion.

Given a weighted particle system $\v w\in[0,1]^n,\ \m X\in\mathbb{B}^{n\times d}$, we denote for $i,j\in\lbrace1,\dots,d\rbrace$ the weighted sample mean by
\begin{align}
\label{eq:sample mean}
\textstyle\bar x_i=\sum_{k=1}^n w_k x_{k,i},\quad \bar x_{i,j}=\sum_{k=1}^n w_k x_{k,i}x_{k,j},
\end{align}
and the weighted sample correlation by
\begin{align}
\label{eq:sample corr}
r_{i,j}=\frac{\bar x_{i,j}-\bar x_i\bar x_j}{\sqrt{\bar x_i(1-\bar x_i)\bar x_j(1-\bar x_j)}}.
\end{align}
For $\varepsilon=0.02$, we define the index set
\begin{equation*}
I\eqdef\lbrace i\in\lbrace1,\dots,d\rbrace \mid \ \bar{x}_i\, \notin\, (\varepsilon,1-\varepsilon)\, \rbrace.
\end{equation*}
which identifies the components which have, according to particle system, a marginal probability close to either boundary of the unit interval.

For the components $i\in I$, we do not consider fitting a logistic regression, but set $L_i=\emptyset$ and draw them independently. Precisely, we set $b_{i,i}=\logit(\bar x_i)$ and $\v b_{i,-i}=\v 0$ which corresponds to logistic model without predictors. Firstly, interactions do not really matter if the marginal probability is excessively small or large. Secondly, these components are prone to cause complete separation in the data or might even be constant.

For the conditional distribution of the remaining components $I^c=\lbrace1,\dots,d\rbrace\setminus I$, we construct parsimonious logistic regressions. For $\delta=0.075$, we define the predictor sets
\begin{equation*}
L_i\eqdef\lbrace j\in\lbrace1,\dots,i-1\rbrace \mid \delta < \abs{r_{i,j}} \rbrace,\quad i\in I^c,
\end{equation*}
which identifies the components with index smaller than $i$ and significant mutual association. Running our examples in Section \ref{sec:simu} with $\delta=0$ show that a saturated logistic regression kernel achieves about the same acceptance rates as a sparse one, while setting $\delta=0.075$ dramatically reduces the computational time we need to calibrate the model.

\paragraph{Fitting the model}
We maximise the log-likelihood function $\ell(\v b)=\ell(\v b\mid\v y, \m Z)$ of a weighted logistic regression model by solving the first order condition $\partial\ell/\partial \v b=\v 0$. We find a numerical solution via Newton-Raphson iterations
\begin{equation}
\label{newton log regression}
-\frac{\partial^2\ell(\v b\up{r})}{\partial \v b \v b\t}(\v
b\up{r+1}-\v b\up{r})
=\frac{\partial\ell(\v b\up{r})}{\partial \v b}, \quad r>0,
\end{equation}
starting at some $\v b\up 0$; see Procedure \ref{algo:fit} for the exact terms. Other updating formulas like Iteratively Reweighted Least Squares or quasi-Newton iterations should work as well.

\floatname{algorithm}{Procedure}
\setcounter{algorithm}{4}
\begin{algorithm}[H]
\caption{Fitting the weighted logistic regressions}
\label{algo:fit}
\begin{algorithmic}
\REQUIRE{$\v w=(w_1,\dots,w_n),\ \m X=(\v x_1,\dots,\v x_n)\t,\ \m B\in\R^{d\times d}$}
\FOR {$i\in I^c$}
  \STATE $\m Z\gets(\m X_{L_i},\v 1),\ \v y\gets\m X_i, \ \v b\up 0\gets\m B_{i,L_i\cup\lbrace i \rbrace}$ 
  \REPEAT
  \STATE \vspace{-2em}\begin{align*}
      p_k&\gets\logit^{-1}(\m Z_k \v b\up{r-1}) &\textbf{ for all } k=1,\dots,n \\
      q_k&\gets p_k(1-p_k) &\textbf{ for all } k=1,\dots,n
    \end{align*}
  \STATE \vspace{-4em}\begin{align*}
      \v b\up{r}\gets &\left(\m Z\t \diag{\v w} \diag{\v q} \m Z+\varepsilon\m I_n\right)^{-1} \times \\
                        &\left(\m Z\t \diag{\v w}\right)
		         \left(\diag{\v q} \m Z\, \v b\up{r-1}+\left(\v y - \v p\right)\right)    
    \end{align*}\vspace{-1.5em}
  \UNTIL {$|b_j\up{r}-b_j\up{r-1}|<10^{-3}$ for all $j$}
  \STATE $\m B_{i,L_i\cup\lbrace i \rbrace}\gets\v b$ \\[0.2em]
\ENDFOR
\RETURN $\m B$
\end{algorithmic}
\end{algorithm}

Sometimes, the Newton-Raphson iterations do not converge because the likelihood function is monotone and thus has no finite maximizer. This problem is caused by data with complete or quasi-complete separation in the sample points \citep{albert_84}. There are several ways to handle this issue.
\begin{enumerate}[(a)]
\item
We just halt the algorithm after a fixed number of iterations and ignore the lack of convergence. Such proceeding, however, might cause uncontrolled numerical problems.
\item
In general, \citet{firth_93} recommends Jeffrey's prior but this option is computationally rather expensive. Instead, we might use a Gaussian prior with variance $1/\varepsilon>0$ which adds a quadratic penalty term $\varepsilon\v b\t \v b$ to the log-likelihood to ensure the target-function is convex.
\item
As we notice that some terms of $\v b_i$ are growing beyond a certain threshold, we move the component $i$ from the set of components with associated logistic regression model $I^c$ to the set of independent components $I$.
\end{enumerate}
In practice, we combine the approaches (c) and (d). In Procedure \ref{algo:fit}, we did not elaborate how to handle non-convergence, but added a penalty term to the log-likelihood, which causes the extra $\varepsilon\m I_n$ in the Newton-Raphson update. Since we solve the update equation via Cholesky factorizations, adding a small term on the diagonal also ensures that the matrix is indeed numerically decomposable.

\paragraph{Starting points}
The Newton-Raphson procedure is known to rapidly converge for starting values $\v b_i\up 0$ not too far from the solution $\v b_i\up *$. In absence of prior information about $\v b_i\up *$, we would naturally start with a vector of zeros and maybe setting $b_{i,i}\up 0=\logit(\bar x_i)$.

In the context of our Sequential Monte Carlo algorithm we can do better than that. Recall that, we constructed a smooth sequence $(\pi_t)_{t=0}^\tau$ of distributions which corresponds to a sequence of proposal distributions $(q_t)_{t=0}^\tau=(q_{\theta_t})_{t=0}^\tau$ which is associated to a sequence $(\theta_t)_{t=0}^\tau$ of parameters.

It significantly speeds up the Newton-Raphson procedure if we choose $\m B_{t-1}$ as starting point for the estimation of $\m B_t$. Indeed, towards the end of the Sequential Monte Carlo algorithm, we fit, for the same precision, a logistic regression in less than four iterations on average when starting at $\m B_{t-1}$, compared to about $13$ iterations on average when starting at zero.

\paragraph{Sampling and evaluating}
In the move step of Sequential Monte Carlo we discussed in Section \ref{sec:move}, we need to sample a proposal state $\v y$ from $q_\theta$ and evaluate the likelihood $q_\theta(y)$ to compute the Metropolis-Hastings ratio \ref{eq:ind mh}. For the logistic regression model $q_{\m B}$, we can do both in one go, see Procedure \ref{algo:sampling}.

\begin{algorithm}[H]
\caption{Sampling from the model}
\label{algo:sampling}
\begin{algorithmic}
\REQUIRE{$\m B$}
\STATE $\v y\gets(0,\dots,0),\ p\gets 1$
\FOR {$i=1\dots,d$}\vspace{0.2em}
  \STATE $q\gets\logit^{-1}(b_{i,i}+\sum_{j\in L_i} b_{i,j}y_j)$ \\[0.2em]
  \STATE \textbf{sample } $\gamma_i\sim\ber_q$ \\[0.2em]
  \STATE $p\gets\begin{cases}
	    p\times q    & \textbf{if }\ \ y_i=1 \\
	    p\times(1-q) & \textbf{if }\ \ y_i=0
            \end{cases}$ \\[0.2em]
\ENDFOR
\RETURN $\v y,\ p$
\end{algorithmic}
\end{algorithm}

\subsection{Why not use a simpler model?}
We briefly justify why we should not use a simpler parametric family for our Sequential Monte Carlo application. Indisputably, the easiest parametric family on $\B^d$ that we can think of is a product model
\begin{equation*}
\textstyle
q_{\v p}(\v \gamma)\eqdef\prod_{i=1}^d \ber_{p_i}(\gamma_i)
\end{equation*}
where $\ber_{p_i(x)}(\gamma)=p_i(x)^\gamma(1-p_i(x))^{1-\gamma}$ denotes a Bernoulli distribution with expected value $p_i(x)\in[0,1]$.

\begin{figure}
\caption{Toy example showing how well the product model $q_{\v p}$ and the logistic regression model $q_{\m B}$ replicate the mass function of a difficult posterior distribution $\pi$.}
\label{fig:toy exa}
\begin{center}
\subfigure[true mass function $\pi(\v\gamma)$]{
\includegraphics[width=0.325\textwidth]{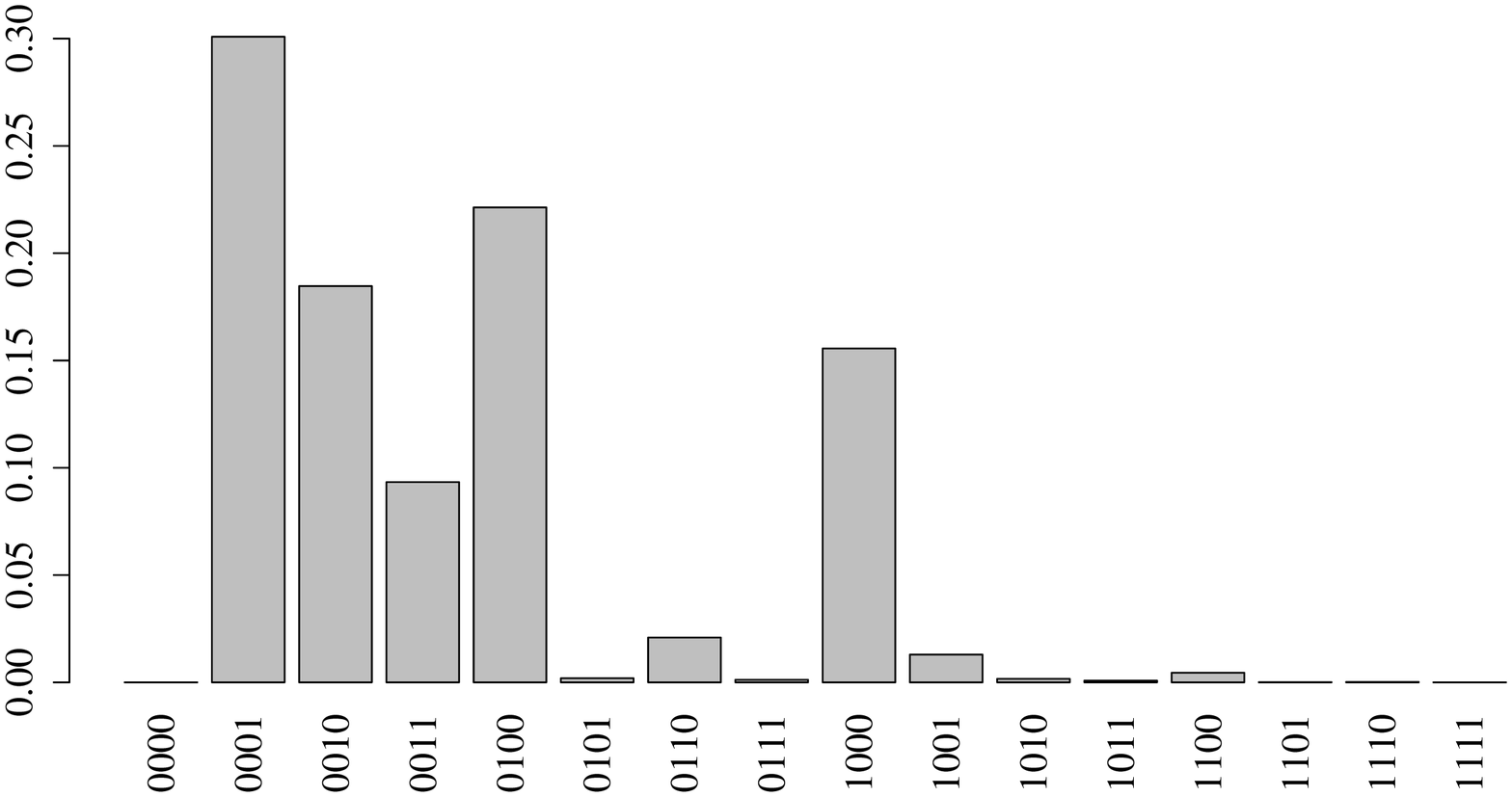}
}\hspace{-4mm}
\subfigure[product model $q_{\v p}(\v\gamma)$]{
\includegraphics[width=0.325\textwidth]{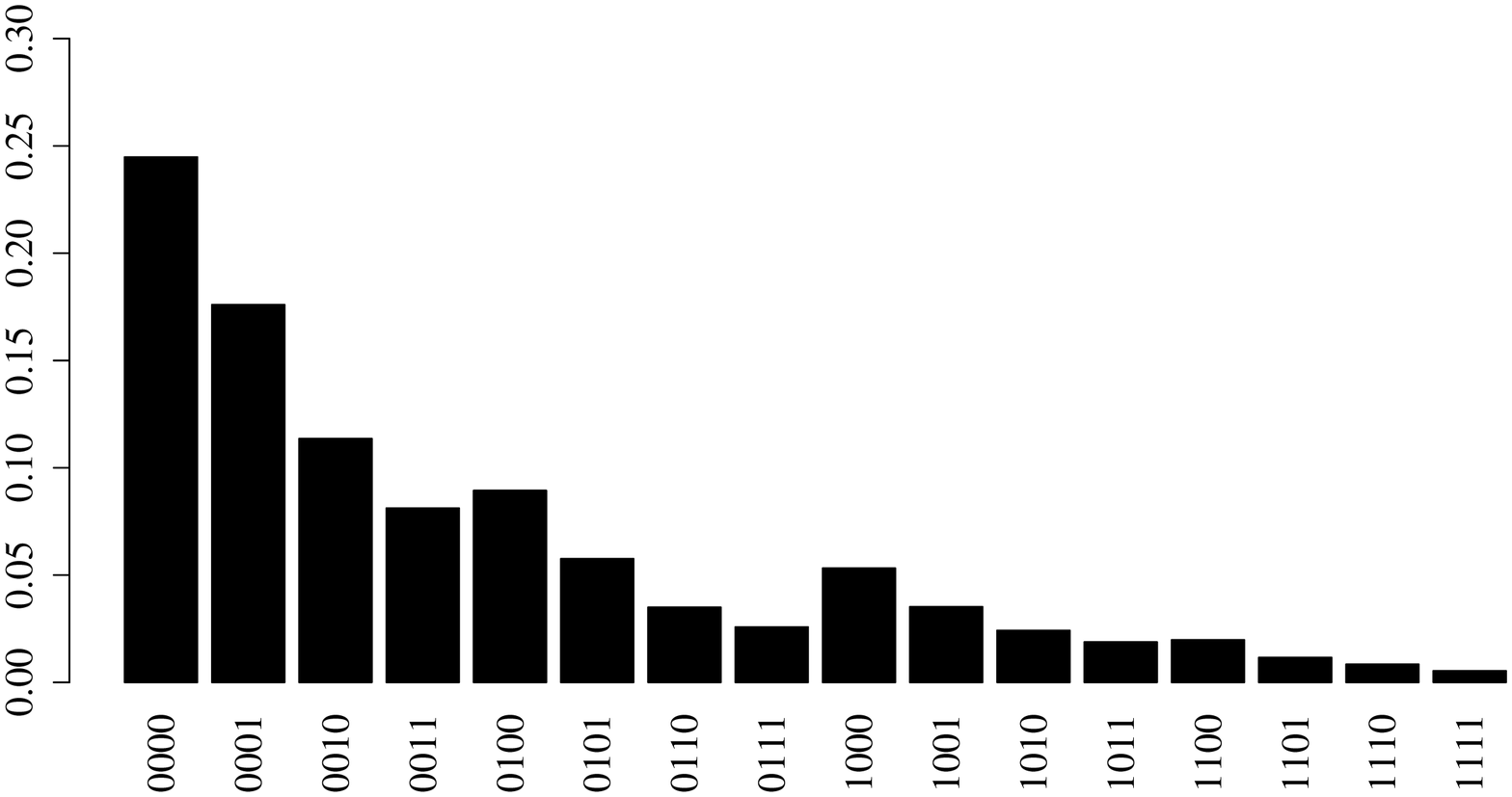}
}\hspace{-4mm}
\subfigure[logistic regression model $q_{\m B}(\v\gamma)$]{
\includegraphics[width=0.325\textwidth]{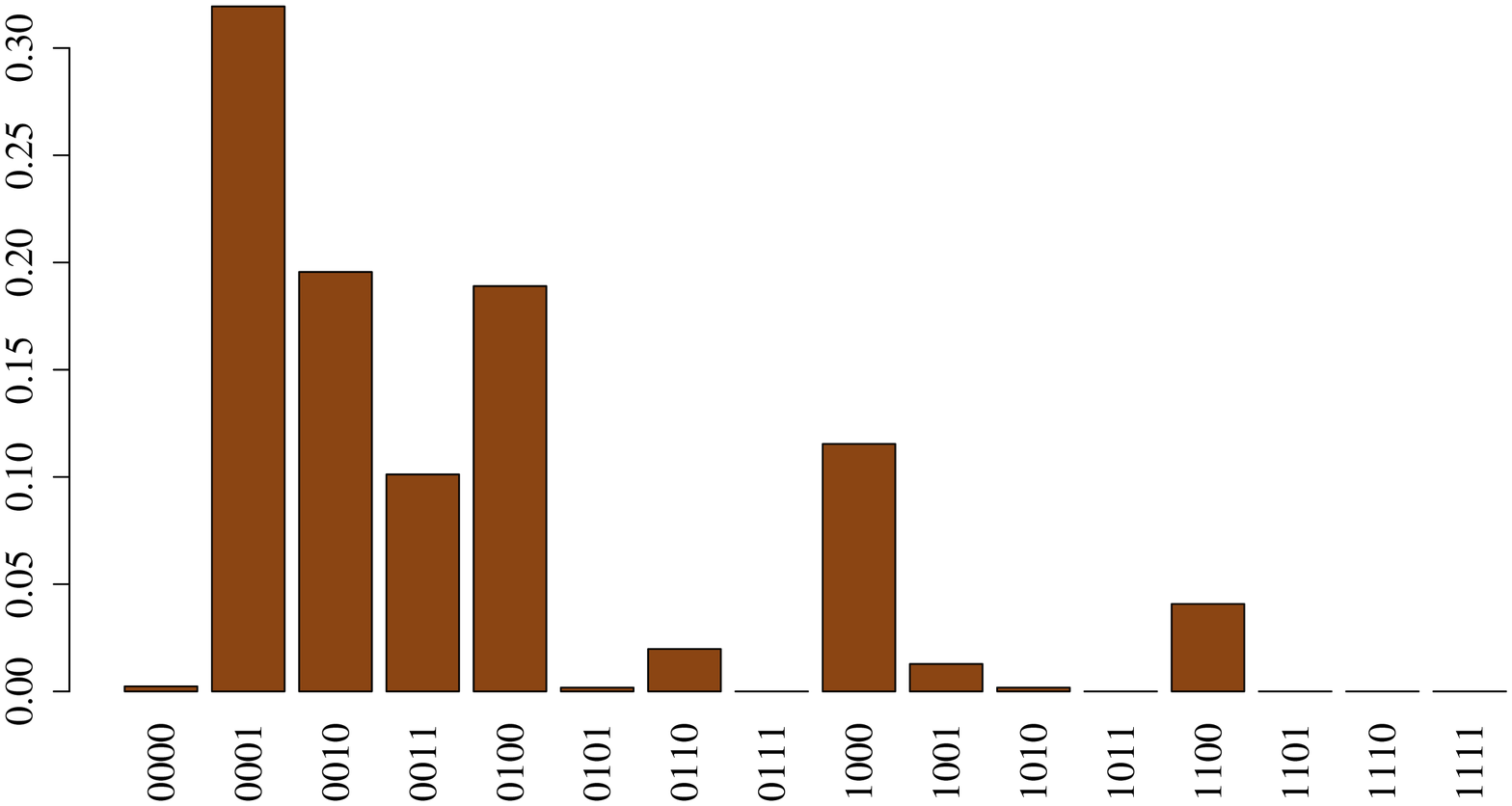}
}
\end{center}
\end{figure}

Let us check the requirement list: the product model is parsimonious with $\mathrm{dim}(\theta)=d$; the maximum likelihood estimator $\theta^*$ is the sample mean $\bar{\v x}=n^{-1}\sum_{k=1}^n x_k$; the decomposition \eqref{eqn:marg dec} holds trivially, which allows us to sample from $q_{\v p}$ and evaluate $q_{\v p}(\v\gamma)$ in $\mathrm{O}(d)$. 

Obviously, however, $q_{\v p}$ does not reproduce any dependencies we might observe in $\m X$. Could we just forget about this last point and use the product model for its simplicity? 

\paragraph{Toy example}
We take a simple linear relation $\v Y=\v V_1+\v V_2$. For $n=100$ and $\mu=10$, we draw normal variates
\begin{equation*}
\v v_1\sim\normal{-\mu}{\m I_n},\ \ \v v_2\sim\normal{\mu}{\m I_n},\ \ \v y=\v v_1+\v v_2
\end{equation*}
where $\v y$ is the vector of observations and
\begin{equation*}
\v z_1,\v z_2\sim\normal{\v v_1}{(\mu^2/4)\,\m I_n},\quad \v z_3,\v z_4\sim\normal{\v v_2}{(\mu^2/4)\,\m I_n}.
\end{equation*}
four columns of covariates.

The posterior distribution $\pi(\v\gamma)=\pi(\v\gamma\mid \v y,\,\m Z)$, using the prior distributions as described in Section \ref{sec:hb}, typically exhibits strong dependencies between its components due to the correlation in the data.

Now we generate pseudo-random data $\m Z$ from $\pi$ and fit both a product model $q_{\v p}$ and a logistic regression model $q_{\m B}$. Looking at the corresponding mass function in Figure \ref{fig:toy exa}, we notice how badly the product model mimics the true posterior. This observation carries over to larger sampling spaces.

\paragraph{Acceptance rates}
A good way to analyse the importance of reproducing the dependencies of $\pi$ is in terms of acceptance rates and particle diversities. As we already remark in Section \ref{sec:move}, the particle diversity naturally diminishes as our particle system approaches a strongly concentrated target distribution $\pi$. However, we want our algorithm to keep up the particle diversity a long as possible to ensure the particle system is well spread out over the entire state space.

\begin{centering}
\begin{figure}
\caption{We compare the use of a product model $q_{\v p}$ to a logistic regression model $q_{\m B}$ as proposal distribution of the Metropolis-Hastings kernel \eqref{eq:ind mh}. We monitor a typical run ($\varrho$ on the x-axis) of our Sequential Monte Carlo algorithm (for the Boston Housing data set described in Section \ref{sec:data sets}) and plot the acceptance rates and particle diversities (on the y-axis).}
\label{fig:ar pd}
\subfigure[acceptance rates]{
\includegraphics[width=0.48\textwidth]{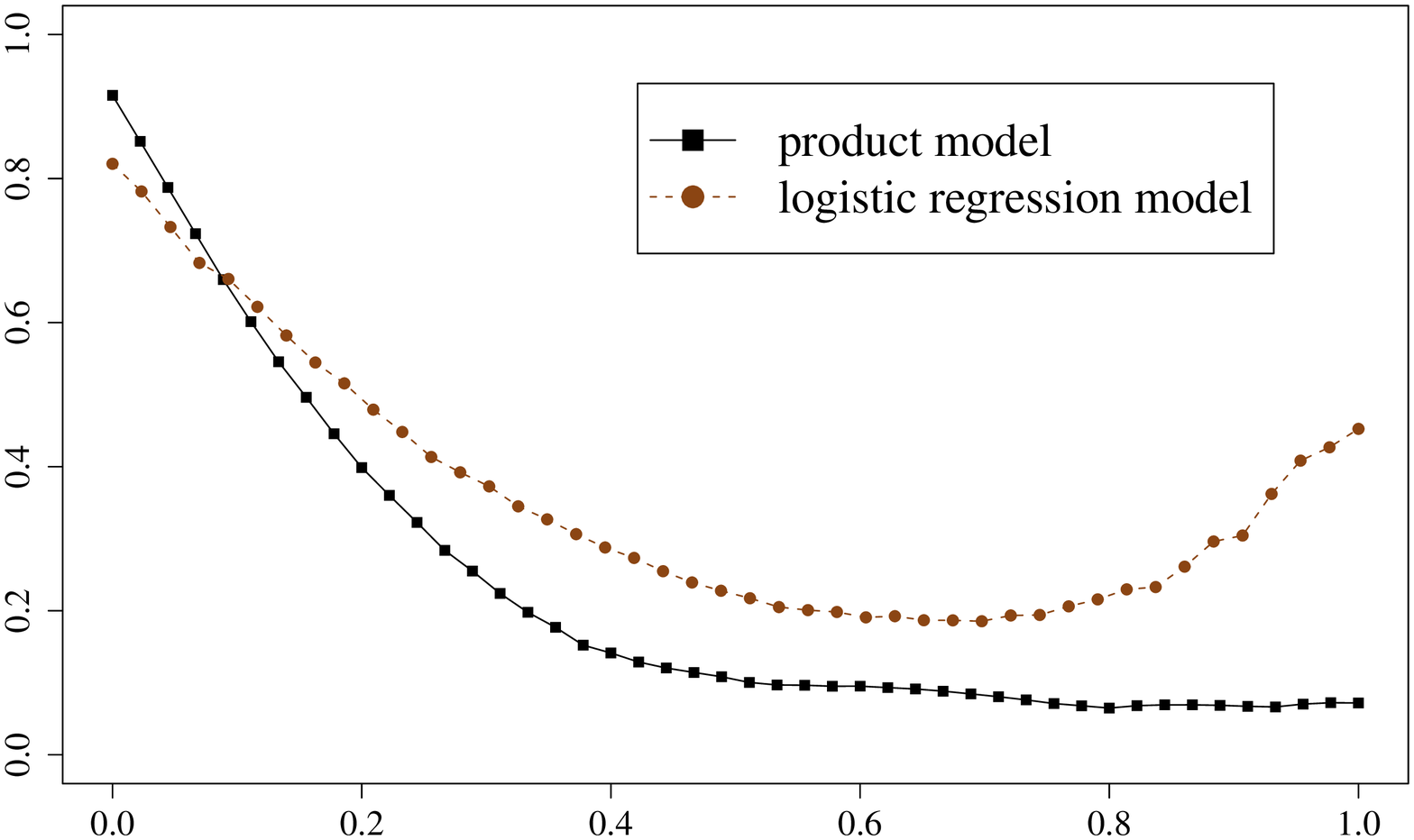}
\label{fig:ar}
}
\subfigure[particle diversities]{
\includegraphics[width=0.48\textwidth]{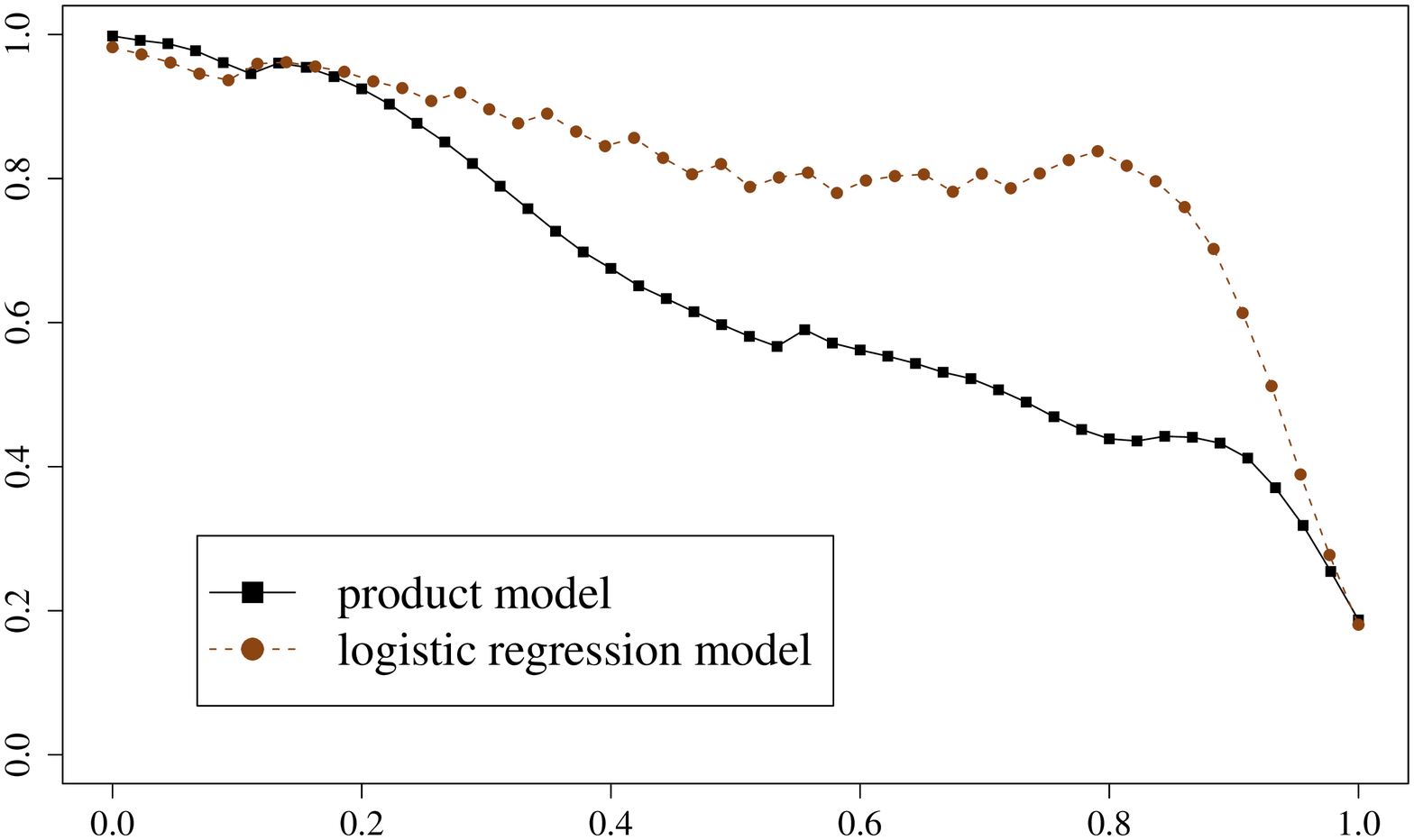}
\label{fig:pd}
}
\end{figure}
\end{centering}

In Figure \ref{fig:ar pd}, we show a comparison (based on the Boston Housing data set explained in Section \ref{sec:data sets}) between two Sequential Monte Carlo algorithms, using a product model $q_{\v p}$ and a logistic regression model $q_{\m B}$ as proposal distribution of the Metropolis-Hastings kernel \eqref{eq:ind mh}.

Clearly, in Figure \ref{fig:ar}, the acceptance rates achieved by the product kernel rapidly decrease and dwell around $5\%$ for the second half of the run. In contrast, the logistic regression kernel always provides acceptance rates greater than $20\%$. As a consequence, in Figure \ref{fig:pd}, the particle diversity sustained by the product kernel decreases at an early stage, while the logistic regression kernel holds it up until the very last steps.

At first sight, it might seem odd that the acceptance rates of the logistic regression kernel increase during the final steps of the algorithm. If we jump ahead, however, and take a look at the results of the Boston Housing problem, see Figure \ref{fig:boston_smc}, we notice that quite a few marginal probabilities of the posterior $\pi$ turn out to be zero, which makes it easier to reproduce the distributions towards the end of the Resample-Move algorithm.

However, if we already decide at an early stage that for some component $\prob{\gamma_i=1}=0$, we fail to ever consider states $\v \gamma\in\B^d$ with $\gamma_i=1$ for the rest of the algorithm. Therefore, the advantage of the logistic regression kernel over the simple product kernel is that we do not completely drop any components from the variable selection problem until the final steps.

\subsection{Review of alternative binary models}
In the following, we review some alternative approaches to modeling multivariate binary data. Unfortunately, we cannot incorporate any of these models in our Sequential Monte Carlo algorithm. Still, it is instructive to understand why alternative strategies fail to provide suitable proposal distributions in the sense of Section \ref{sec:model prop}. For a more detailed review of parametric families suitable for adaptive Monte Carlo algorithms on binary spaces, see \citet{schaefer2011parametric}.

\paragraph{Quadratic multi-linear models}
For coefficients $a\in\R^{2^d}$, we can write any mass function on $\B^d$ as
\begin{equation*}
\textstyle
\pi(\v\gamma)=\sum_{S\subseteq \lbrace1,\dots,d\rbrace} a_S\prod_{i\in S}\gamma_i.
\end{equation*}
It is tempting to construct a $d(d+1)/2$ parameter model
\begin{equation*}
q_{\mu,\m A}(\v\gamma)\eqdef\mu + \v\gamma\t\m A\v\gamma
\end{equation*}
by removing interaction terms of order higher than two. As \cite{bahadur61representation} points out, the main problem of any additive approach is the fact that a truncated model might not be non-negative and thus not define a probability distribution.

Although the linear structure allows to derive explicit and recursive formulae for the marginal and conditional distributions, we hardly ever find a useful application for the additive model. As other authors \citep{park1996simple,emrich1991method} remark, additive representations like the much-cited \cite{bahadur61representation} expansion are quite instructive but, unfortunately, impractical.

\paragraph{Quadratic exponential models}
For coefficients $a\in\R^{2^d}$, we can write any mass function on $\B^d$ as
\begin{equation*}
\textstyle
\pi(\v\gamma)=\exp\left(\sum_{S\subseteq \lbrace1,\dots,d\rbrace} a_S\prod_{i\in S}\gamma_i.\right)
\end{equation*}
Removing higher order interaction terms, we can construct a $d(d+1)/2$ parameter model
\begin{equation}
\label{eqn:eqm}
q_{\mu,\m A}(\v\gamma)\eqdef\mu\exp(\v\gamma\t\m A\v\gamma),
\end{equation}
where $\m A$ is a symmetric matrix. Quadratic exponential models are a well defined class of distributions, but there is no simple recursive structure for their marginal distributions. Hence, we cannot compute the factorization \eqref{eqn:marg dec} we need to sample from $q_{\m A}$.

\cite{cox1994note} propose an approximation to the marginal distributions by expressions of the form \eqref{eqn:eqm}, omitting higher order terms in a Taylor expansion. If we write the parameter $\m A$ as
$$
\m A=\begin{pmatrix}
	\m A' & \v b\t \\
	\v b  & c
  \end{pmatrix},
$$
the parameter of the marginal distribution $q_{\m A_{1:d-1}}(\v\gamma_{\,1:d-1})$ is approximately given by
\begin{equation*}
\textstyle
\m A_{1:d-1}\approx
\m A'
+\left(1+\tanh(\frac{c}{2})\right)\diag{\v b}
+\frac{1}{2}\,\mathrm{sech}^2(\frac{c}{2})\v b \v b\t,
\end{equation*}
and the normalizing constant is $\mu_{1:d-1}=\mu(1+\exp(c))$. We can recursively compute approximations to all marginal distributions $q_{\m A_{1:d-1}},\dots,q_{\m A_{1:1}}$ and derive logistic forms
\begin{equation*}
\logit(\prob{\gamma_i=1\mid\v \gamma_{1:i-1}})=
\log \frac{q_{A_{1:i}}(\gamma_i=1,\v \gamma_{1:i-1})}{q_{A_{1:i}}(\gamma_i=0,\v \gamma_{1:i-1})},
\end{equation*}
which takes us back to \eqref{eqn:lb}. However, there is no reason to fit a quadratic exponential model and compute the approximate logistic model if we can directly fit the logistic conditionals model in the same time.

\paragraph{Latent variable models}
Let $\varphi_\theta$ be a parametric family on $\mathbb X$ and $\tau\colon\mathbb X \to \B^d$ a mapping into the binary state space. We can sample from a latent variable model
\begin{equation*}
\textstyle
q_{\theta}(\v\gamma)\eqdef\int_{\tau^{-1}(\v\gamma)}\,\varphi_\theta(v)\,d\v v
\end{equation*}
by setting $\v y=\tau(\v v)$ for a draw $\v v\sim \varphi_\theta$ from the latent parametric family.

Non-normal parametric families with $d(d-1)/2$ dependence parameters seem to either have a very limited dependence structure or unfavourable properties \citep{joe1996families}. Therefore, the multivariate normal
\begin{align*}
\varphi_{(\v \mu,\m \Sigma)}(\v v)&=(2\pi)^{-d/2}\abs{\m \Sigma}^{-1/2}e^{-1/2(v-\v \mu)\t\m \Sigma^{-1}(v-\v \mu)}, \\
\tau(\v v)&=(\ind_{(\infty,0]}(v_1),\dots,\ind_{(\infty,0]}(v_d)),
\end{align*}
appears to be the natural and almost the only option for $p_\theta$. This kind of model has been discussed repeatedly in the literature \citep{emrich1991method, leisch1998generation, cox2002some}.

The first and second order marginal probabilities of the model $q_{(\v \mu,\m \Sigma)}$ are given by $\varPhi_1(\mu_i)$ and $\varPhi_2(\mu_i,\mu_j;\sigma_{i,j})$, respectively, where $\varPhi_1(v_i)$ and $\varPhi_2(v_i,v_j;\sigma_{i,j})$ denote the cumulative distribution functions of the univariate and bivariate normal distributions with zero mean, unit variance and correlation $\sigma_{i,j}\in[-1,1]$.

We can fit the model $q_{(\v\mu,\m \Sigma)}$ to a particle system $(\v w,\m X)$ by matching the moment, that is adjusting $\v \mu$ and $\m \Sigma$ such that
\begin{equation*}
\varPhi_1(\mu_i)=\bar x_i,\quad \varPhi_1(\mu_i,\mu_j;\sigma_{i,j})=r_{i,j}
\end{equation*}
with $\bar x_i$ and $r_{i,j}$ as defined in \eqref{eq:sample mean} and \eqref{eq:sample corr}. However, the locally constructed correlation matrix $\m \Sigma$ might not be positive definite. Still, we can obtain a feasible parameter replacing $\m \Sigma$ by $\m \Sigma^*=(\m \Sigma+\abs{\lambda}\m I)/(1+\abs{\lambda})$, where $\lambda$ is smaller than all eigenvalues of the locally adjusted matrix $\m \Sigma$.

The main drawback of latent variable approaches is the fact that that point-wise evaluation of the probability mass function $q_\theta(\v y)$ is computationally feasible only in special cases. Hence, we cannot use this class of models in a Sequential Monte Carlo context.

\paragraph{Archimedean copula models}
The potentials and pitfalls of applying copula theory, which is well developed for bivariate, continuous random variables, to multivariate discrete distribution is discussed in \citet{genest2007primer}. There have been earlier attempts to sample binary vectors via copulae: \citet{lee1993generating} describes how to construct an Archi\-medean copula, more precisely the Frank family \citep[p.119]{nelsen2006introduction}, for sampling multivariate binary data. Unfortunately, this approach is limited to very low dimensions.

\paragraph{Multivariate reduction models}
Several approaches to generating multivariate binary data are based on a representation of the components $\v\gamma_i$ as functions of sums of independent variables \citep{park1996simple, lunn1998note, oman2001modelling}. These techniques are limited to certain patterns of non-negative correlation, and do, therefore, not yield suitable proposal distributions in a Sequential Monte Carlo application. We mention them for the sake of completeness.

%

%


%

%

\section{Numerical experiments}
\label{sec:simu}
In this section we compare our Sequential Monte Carlo algorithm to standard Markov chain methods based on local moves as introduced in Section \ref{sec:mcmc}. These are standard algorithms and widely used. There are other recent approaches like Bayesian Adaptive Sampling \citep{clyde2011bayesian} or Evolutionary Stochastic Search \citep{bottolo2010ess} which also aim at overcoming the difficulties of multi-modal binary distributions. However, a thorough and just comparison of our Sequential Monte Carlo approach to other advanced methods needs careful consideration and is beyond the scope of this paper.

For testing, we created variable selection problems with high dependencies between the covariates which yield particularly challenging, multi-modal posterior mass functions. The problems are build from freely available datasets by adding logarithms, polynomials and interaction terms. The Markov chain Monte Carlo methods presented in Section \ref{sec:mcmc} tend to fail on these problems due to the very strong multi-modality of the posterior distribution while the Sequential Monte Carlo approach we advocate in Section \ref{sec:smc} yields very reliable results.

Note, however, that using Sequential Monte Carlo we do not get something for nothing. Firstly, the implementation of our algorithm including the logistic conditionals model introduced in Section \ref{sec:logistic model} is quite involved compared to standard Markov chain algorithms. Secondly, simple Markov chain methods are faster than our algorithm while producing results of the same accuracy if the components of the target distribution are nearly independent.

\subsection{Construction of the data sets}
\label{sec:data sets}
We briefly describe the variable selection problems composed for our numerical experiments.

\paragraph{Boston Housing}
The first example is based on the Boston Housing data set, originally treated by \cite{harrison1978hedonic}, which is freely available at the \href{http://lib.stat.cmu.edu/datasets/boston\_corrected.txt}{StatLib} data archive. The data set provides covariates ranging from the nitrogen oxide concentration to the per capita crime rate to explain the median prices of owner-occupied homes, see Table \ref{tab:bostonsum}. The data has yet been treated by several authors, mainly because it provides a rich mixture of continuous and discrete variables, resulting in an interesting variable selection problem.

Specifically, we aim at explaining the logarithm of the corrected median values of owner-occupied housing. We enhance the $13$ columns of the original data set by adding first order interactions between all covariates. Further, we add a constant column and a squared version of each covariate (except for \textsc{chas} since it is binary).

This gives us a model choice problem with $104$ possible predictors and $506$ observations. We use a hierarchical Bayesian approach, with priors as explained in the above Section \ref{sec:hb}, to construct a posterior distribution $\pi$. By construction, there are strong dependencies between the possible predictors which leads to a rather complex, multi-modal posterior distribution.

\begin{table}[hb]
\vspace{-5mm}
\caption{Boston Housing data summary.}
\label{tab:bostonsum}
\begin{center}
\begin{tabular}{ll}
short name & explanation\\
\noalign{\smallskip}\hline\noalign{\smallskip}
\textsc{crim}    & per capita crime\\
\textsc{zn}      & proportions of residential land zoned\\
& for lots over $2323\ \mathrm{m}^2$\\
\textsc{indus}   & proportions of non-retail business acres\\
\textsc{chas}    & tract borders Charles River (binary) \\
\textsc{nox}     & nitric oxides concentration (parts per $10^7$)\\
\textsc{rm}      & average numbers of rooms per dwelling\\
\textsc{age}     & proportions of owner-occupied units\\
& built prior to 1940\\
\textsc{dis}     & weighted distances to five Boston\\
& employment centres\\
\textsc{rad}     & accessibility to radial highways\\
\textsc{tax}     & full-value property-tax rate per USD $10^4$\\
\textsc{ptratio} & pupil-teacher ratios\\
\textsc{b}       & $(\mathrm{Bk} - 0.63)^2$ where $\mathrm{Bk}$ is the proportion\\
& of the black population\\
\textsc{lstat}   & percentage of lower status population\\
\end{tabular}
\end{center}
\vspace{-10mm}
\end{table}

\paragraph{Concrete Compressive Strength}
The second example is constructed from a less known data set, originally treated by \cite{yeh1998modeling}, which is freely available at the \href{http://archive.ics.uci.edu/ml/machine-learning-databases/concrete/compressive/Concrete\_Data.xls}{UCI Machine Learning Repository}. The data provides information about components of concrete to explain its compressive strength. The compressive strength appears to be a highly non-linear function of age and ingredients.

In order to explain the compressive strength, we take the $8$ covariates of the original data set and add the logarithms of some covariates (indicated by the prefix \textsc{lg}), see Table \ref{tab:concretesum}. Further, we add interactions between all $13$ covariates of the augmented data set and a constant column.

This gives us a model choice problem with $79$ possible predictors and $1030$ observations. We use a hierarchical Bayesian approach, with priors as explained in the above Section \ref{sec:hb}, to construct a posterior distribution $\pi$.

\begin{table}[hb]
\vspace{-5mm}
\caption{Concrete Compressive Strength data summary. Components are measured as kg/$\mathrm{m}^3$.}
\label{tab:concretesum}
\begin{center}
\begin{tabular}{ll}
short name & explanation\\
\noalign{\smallskip}\hline\noalign{\smallskip}
\textsc{c}, \textsc{lg\_c}      & cement\\
\textsc{blast}                  & blast furnace slag\\
\textsc{fash}                   & fly ash\\
\textsc{w}, \textsc{lg\_w}      & water\\
\textsc{plast}                  & superplasticizer\\
\textsc{ca}, \textsc{lg\_ca}    & coarse aggregate\\
\textsc{fa}, \textsc{lg\_fa}    & fine aggregate\\
\textsc{age}, \textsc{lg\_age}  & age in days\\
\end{tabular}
\end{center}
\vspace{-10mm}
\end{table}

\paragraph{Protein activity data}
The third example has originally been analyzed by \cite{clyde1998protein}. Later, \cite{clyde2011bayesian} used it as a challenging example problem in variable selection and included the raw data in the \textsc{R}-package \href{http://cran.r-project.org/web/packages/BAS/index.html}{BAS} available at \href{http://cran.r-project.org/web/packages/BAS/index.html}{CRAN} which implements the Bayesian Adaptive Sampling algorithm.

In order to explain the protein activity (\textsc{prot.act1}), we first convert the factors \textsc{buf}, \textsc{ra} and \textsc{det} into a factor model. We enhance the $14$ columns of this data set by adding first order interactions between all covariates and a constant column. For details on the raw data see Table \ref{tab:proteinsum}.

Note that some columns turn out to be constant zeros such that we obtain a model choice problem with $88$ possible predictors and $96$ observations. For reasons of consistency, we choose the priors explained in the above Section \ref{sec:hb} instead of the original $g$-prior used in \cite{clyde2011bayesian}.

\begin{table}[hb]
\vspace{-5mm}
\caption{Protein activity data summary.}
\label{tab:proteinsum}
\begin{center}
\begin{tabular}{ll}
short name & explanation \\
\noalign{\smallskip}\hline\noalign{\smallskip}
\textsc{det}   & detergent \\
\textsc{buf}   & pH buffer \\
\textsc{NaCl}  & salt \\
\textsc{con}   & protein concentration \\
\textsc{ra}    & reducing agent  \\
\textsc{MgCl2} & magnesium chloride \\
\textsc{temp}  & temperature \\
\end{tabular}
\end{center}
\vspace{-10mm}
\end{table}

\subsection{Main effect restrictions}
In some statistical applications we might want to only include the interactions if the corresponding main effects are present in the model. These constraints are easy to incorporate if needed but render the sampling problem even more challenging since the constrained support makes the state space exploration more difficult.

Let $d$ denote the number of main effects and $\gamma_{i,j}$ the interaction of the main effects $\gamma_i$ and $\gamma_j$ for all $i,j=1,\dots,d$. For the variable selection problem on $\B^{d(d+1)/2}$, we impose the prior constraints on the feasible interactions
\begin{equation}
\label{sec:me prior}
\pi(\v\gamma\mid\m Z)\propto\ind_{\{\v\gamma\in\B^{d(d+1)/2}\mid \gamma_{ij}\leq\gamma_i\gamma_j \text{ for all } i,j=1,\dots,d\}}(\v\gamma).
\end{equation}
While the Markov chain Monte Carlo algorithms can proceed as before, we need to slightly modify our Sequential Monte Carlo approach. When sampling from a restricted distribution, we initialise the particle system with an iid sample from the prior \eqref{sec:me prior} instead of the uniform distribution. In the sequel, we report for each dataset a comparison with and without the main effect restrictions.

\subsection{How to compare to Markov chain Monte Carlo}
We do not think it is reasonable to compare two completely different algorithms in terms of pure computational time. We cannot guarantee that our implementations are optimal nor that the time measurements can exactly be reproduced in other computing environments.

We suppose that the number of evaluations of the target function $\pi$ is more of a fair stopping criterion, since it shows how well the algorithms exploit the information obtained from $\pi$. Precisely, we parameterise the Sequential Monte Carlo algorithm to not exceed a fixed number $\nu$ of evaluations and stop the Markov chains when $\nu$ evaluations have been performed.

\paragraph{Assets and drawbacks}
The Sequential Monte Carlo and the Markov chain Monte Carlo algorithms both have extensions and numerical speed-ups which make it hard to settle on a fair comparison.

Advocates of  Markov chain methods might criticise that the number of target evaluations is a criterion biased towards the Sequential Monte Carlo approach, for there are updating schemes which allow for faster computation of the Cholesky decomposition \eqref{eq:chol} given the decomposition of a neighbouring model, see \citet[chaps. 8,10]{dongarra1979linpack}. Thus, Markov chains which propose to change one component in each step can evaluate $\pi$ with less effort and perform more evaluations of $\pi$ in the same time.

On the other hand, however, the Sequential Monte Carlo algorithm can be parallelised in the sense that we can, on suitable hardware, run many evaluations of $\pi$ in parallel during the move step, see Procedure \ref{algo:move}. No analogue speed-up can be performed in the context of Markov chains. We have processed variable selection problems from genetics with about a thousand covariates within a few hours running a parallelised version of the algorithm on a $64$-CPU cluster. A detailed report is going to be published as supplementary material.

Further, Sequential Monte Carlo methods are more suitable than Markov chain Monte Carlo to approximate the evidence, that is the normalization constant of the posterior distribution. We can exploit this property to compare, for instance, regression models with different monotonic link functions.

\paragraph{Parameters}
We run our Sequential Monte Carlo (SMC) algorithm with $n=1.5\times10^4$ particles and a target effective sample size $\eta=0.9$, as explained in Section \ref{sec:smc}. For these parameters, the Sequential Monte Carlo algorithm needs less than $\nu=2.5\times10^6$ evaluations of $\pi$ on all examples problems.

We compare our algorithm to both the Adaptive Markov chain Monte Carlo \citep[AMCMC]{nott2005adaptive} and the standard metropolised Gibbs \citep[MCMC]{liu1996peskun} described in Section \ref{sec:mcmc}. For the MCMC, we draw the number of bits to be flipped from a truncated geometric distribution with mean $k^*=2$ as proposed in Section \ref{sec:blockup}. However, as stated earlier, we could not observe a significant effect of changes in the block updating schemes on the quality of the Monte Carlo estimate.

For the AMCMC, we use $\delta=0.01$ and $\lambda=0.01$, following the recommendations of \citet{nott2005adaptive}. We update the estimates $\psi$ and $\m W$ every $2\times10^5$ iterations of chain. Before we start adapting, we generate $2.5\times10^5$ iterations with a metropolised Gibbs kernel (after a discarded burn-in of $2.5\times10^4$ iterations).

\subsection{Implementation}
The numerical work was completely done in \href{http://www.python.org/}{Python 2.6} using \href{http://www.scipy.org}{SciPy} packages and run on a cluster with $1.86$GHz processors. Scientific work in applied fields is often more accessible to the reader if the source code which generated numerical evidence is released along with the publication. The complete, documented sources used in this work can be found at \url{http://code.google.com/p/smcdss}.

We also provide instructions on how to install and run our project. The program can process data sets in standard \href{http://en.wikipedia.org/wiki/Comma-separated_values}{csv}-format and generate \href{http://www.r-project.org/}{R} scripts for graphical visualisation of the results. The released version was tested to run on both Windows and Linux machines.

\subsection{Results and discussion}
\label{sec:results}
We run each algorithm $200$ times and each time we obtain for all covariates a Monte Carlo estimate of the marginal probability of inclusion in the normal linear model. We visualize the variation of the estimator by box-plots that show how much the Monte Carlo estimates have varied throughout the $200$ runs (Figures \ref{fig:boston} and \ref{fig:concrete}). Here, the white boxes contain $80\%$ of the Monte Carlo results, while the black boxes show the extent of the $20\%$ outliers. For better readability, we add a coloured bar up to the smallest estimate we obtained in the test runs; otherwise components with a small variation are hard to see.

The vertical line in the white box indicates the median of the Monte Carlo estimates. The median of the Sequential Monte Carlo runs correspond very precisely to the results we obtained by running a Markov chain Monte Carlo algorithm for a few days. Unquestionably, the Sequential Monte Carlo algorithm is extremely robust; for $200$ test runs and for both data sets, the algorithm did not produce a single major outlier in any of the components.

This not true for either of the Markov chain algorithms. The size of white boxes indicate that adaptive Markov chain Monte Carlo works quite better than the standard Markov chain procedure. However, even the adaptive Markov chain method is rather vulnerable to generating outliers. The large black boxes indicate that, for some starting points of the chain, the estimates of some marginal probabilities might be completely wrong.

The outliers, that is the black boxes, in the MCMC and the AMCMC plots are strikingly similar. The adaptive and the standard Markov chains apparently both fall into the same trap, which in turn confirms the intuition that adaption makes a method faster but not more robust against outliers. An adaptive local method is still a local method and does not yield reliable estimates for difficult binary sampling problems. Figure \ref{fig:me protein} suggests that in constrained spaces adaption is difficult and might even have contra-productive effects.

In Tables \ref{fig:boston} to \ref{fig:me protein}, we gather some key performance indicators, each averaged over the $200$ runs of the respective algorithms. Note that the time needed to perform $2.5\times10^6$ evaluations of $\pi$ is a little less than the running time of the standard Markov chain. Thus, even in terms of computational time, the adaptive Markov chain can hardly compete with our Sequential Monte Carlo method, even if evaluations of $\pi$ were at no cost.

%


%

\def\histwidth{0.31\textwidth}
\begin{center}
\begin{figure*}
\vspace{-2mm}
\caption{Boston Housing data set. For details see Section \ref{sec:results}.}
\label{fig:boston}
\subfigure[SMC $\sim1.4\times10^6$eval'ns of $\pi$]{
\includegraphics[angle=-90,width=\histwidth]{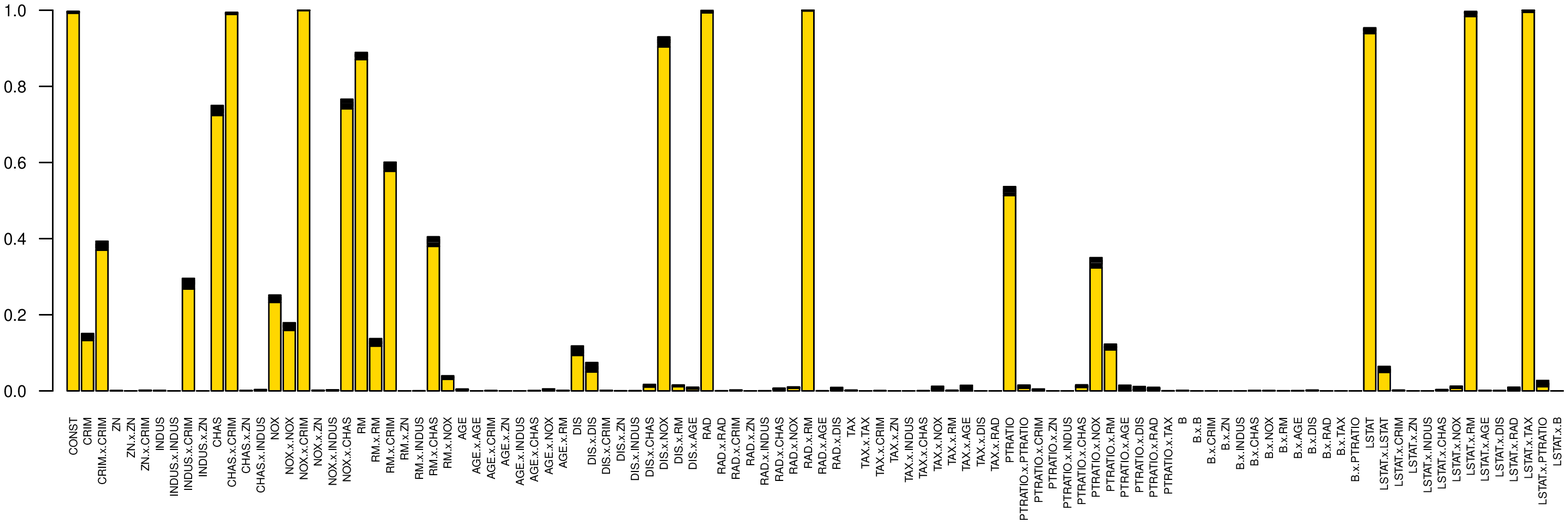}
\label{fig:boston_smc}
}
\subfigure[AMCMC $2.5\times10^6$eval'ns of $\pi$]{
\includegraphics[angle=-90,width=\histwidth]{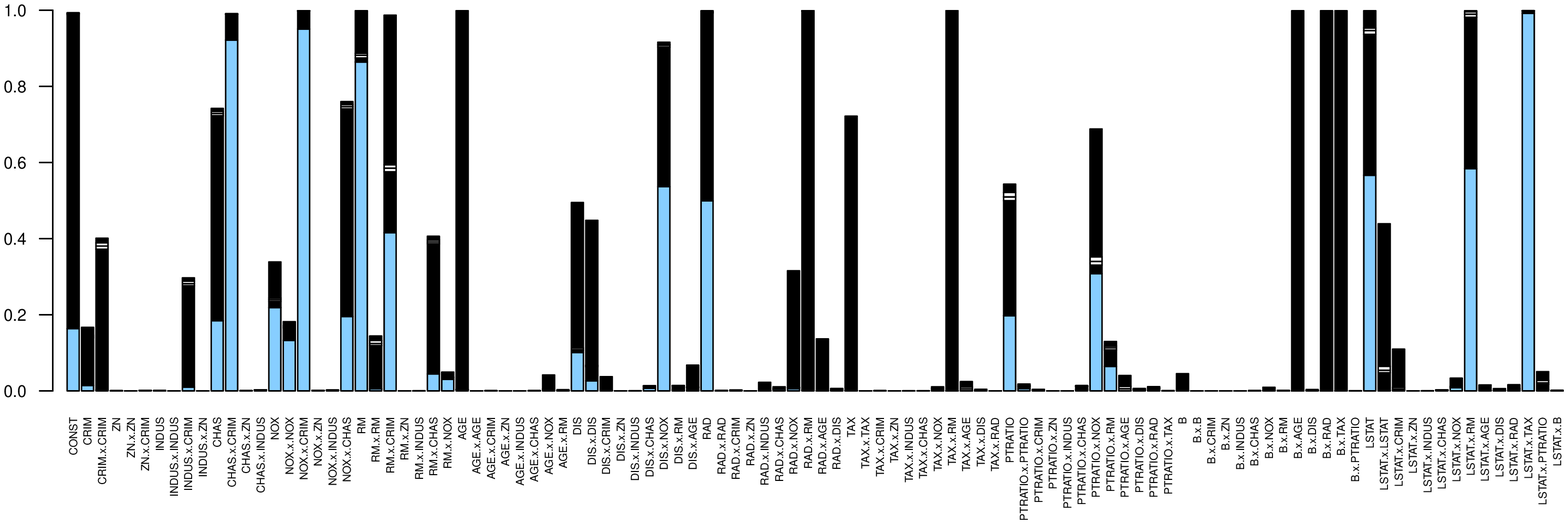}
\label{fig:boston_amcmc}
}
\subfigure[MCMC $2.5\times10^6$eval'ns of $\pi$]{
\includegraphics[angle=-90,width=\histwidth]{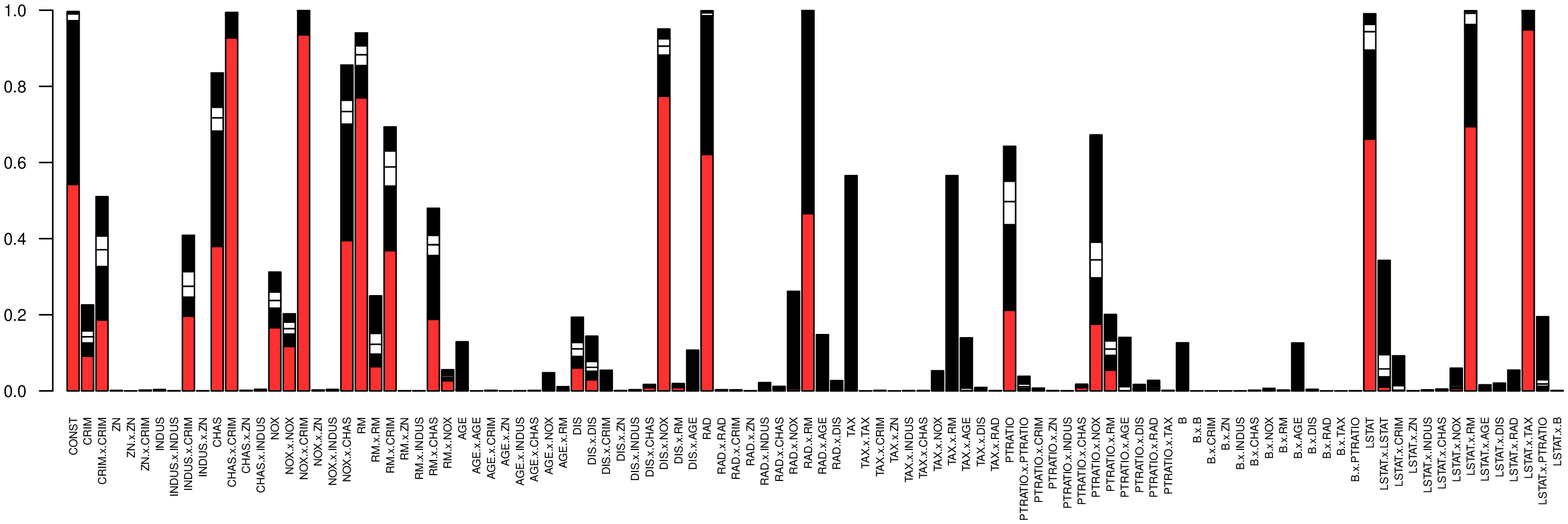}
\label{fig:boston_mcmc}
}
\vspace{1mm}
\begin{flushleft}
\small{\bfseries Table.} Boston Housing data set. Averaged key indicators complementary to Figure \ref{fig:boston}.
\end{flushleft}
\begin{center}
\begin{tabular}[hb]{l|rrr}
& Sequential MC & Adaptive MCMC & Standard MCMC \\
\noalign{\smallskip}\hline\noalign{\smallskip}
computational time & $0:36:59$ h & $4:50:52$ h & $0:38:06$ h \\
evaluations of $\pi$ & $1.36\times10^6$ & $2.50\times10^6$ & $2.50\times10^6$ \\
average acceptance rate & $36.4\%$ &  $29.1\%$  &  $0.81\%$ \\
length $t$ of the chain $\v x_t$ & & $7.52 \times 10^7$ & $2.50\times10^6$ \\
moves $\v x_t\neq \v x_{t-1}$ & & $7.28 \times 10^5$ & $2.07\times10^4$ \\
\end{tabular}
\end{center}
\end{figure*}
\end{center}

\begin{center}
\begin{figure*}
\vspace{-2mm}
\caption{Boston Housing data set with main effect restrictions. For details see Section \ref{sec:results}.}
\label{fig:me boston}
\subfigure[SMC $\sim1.2\times10^6$eval'ns of $\pi$]{
\includegraphics[angle=-90,width=\histwidth]{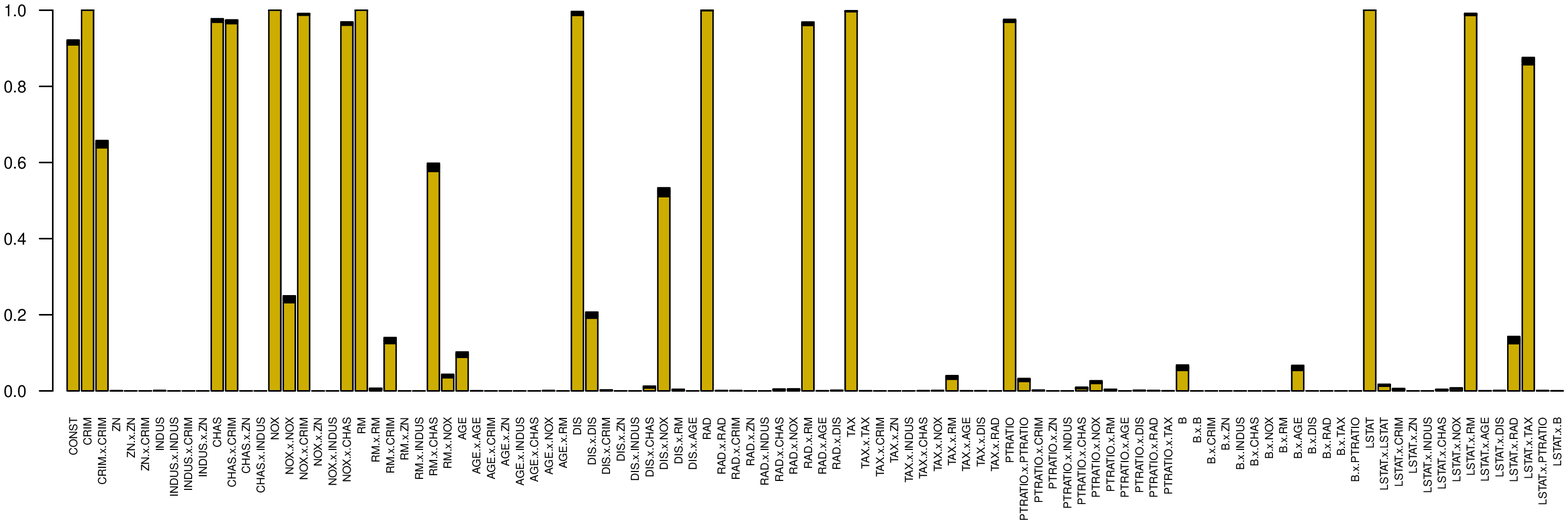}
\label{fig:me boston_smc}
}
\subfigure[AMCMC $2.5\times10^6$eval'ns of $\pi$]{
\includegraphics[angle=-90,width=\histwidth]{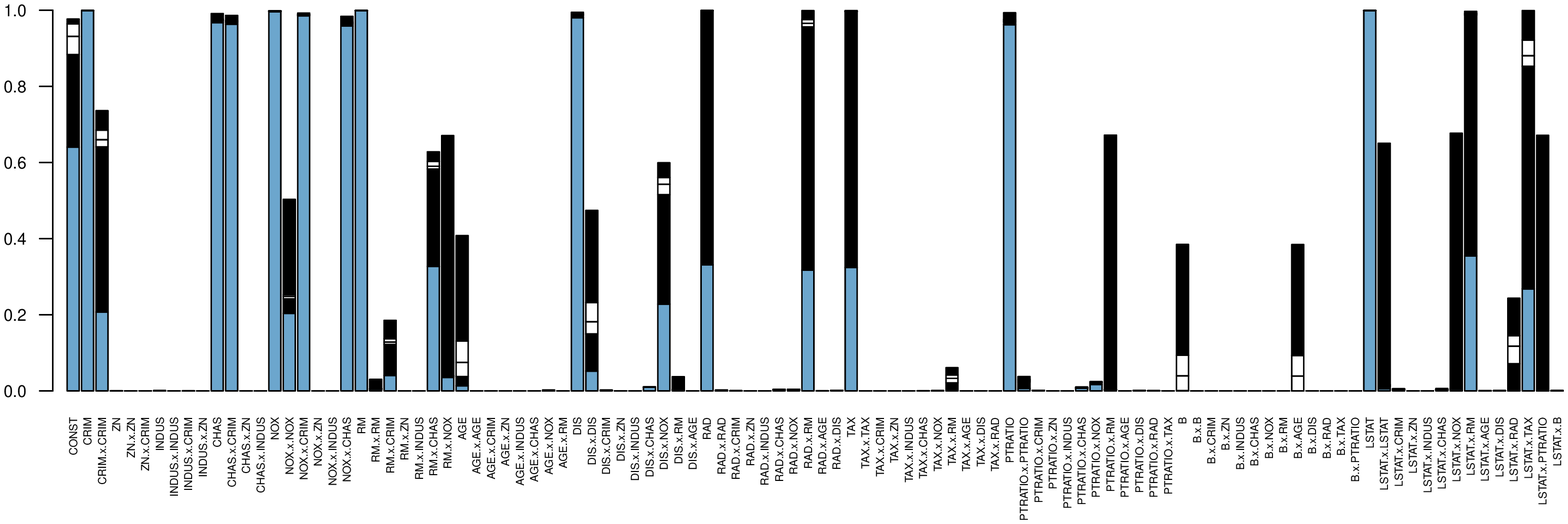}
\label{fig:me boston_amcmc}
}
\subfigure[MCMC $2.5\times10^6$eval'ns of $\pi$]{
\includegraphics[angle=-90,width=\histwidth]{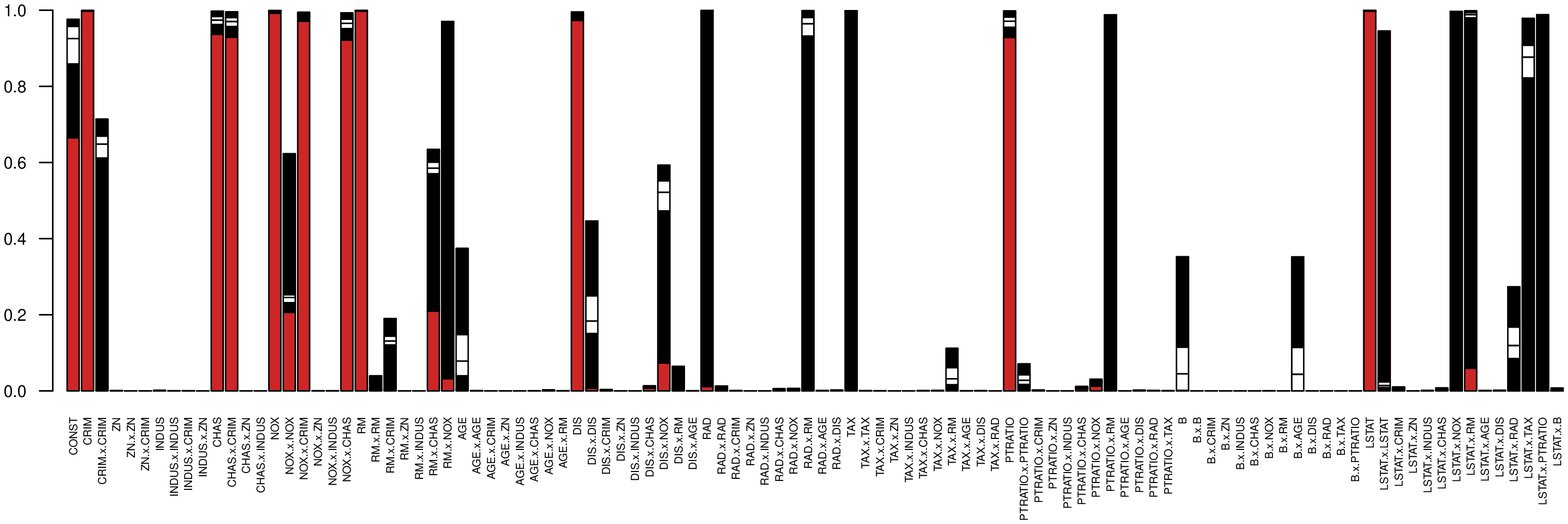}
\label{fig:me boston_mcmc}
}
\vspace{1mm}
\begin{flushleft}
\small{\bfseries Table.} Boston Housing data set with main effect restrictions. Averaged key indicators complementary to Figure \ref{fig:me boston}.
\end{flushleft}
\begin{center}
\begin{tabular}[hb]{l|rrr}
& Sequential MC & Adaptive MCMC & Standard MCMC \\
\noalign{\smallskip}\hline\noalign{\smallskip}
computational time & $0:18:05$ h & $4:33:20$ h & $0:14:13$ h \\
evaluations of $\pi$ & $1.15\times10^6$ & $2.50\times10^6$ & $2.50\times10^6$ \\
average acceptance rate & $20.79\%$ &  $45.4\%$  &  $1.20\%$ \\
length $t$ of the chain $\v x_t$ & & $8.01 \times 10^7$ & $2.50\times10^6$ \\
moves $\v x_t\neq \v x_{t-1}$ & & $1.13 \times 10^6$ & $2.96\times10^4$ \\
\end{tabular}
\end{center}
\end{figure*}
\end{center}

%


%

\begin{center}
\begin{figure*}
\vspace{-2mm}
\caption{Concrete Compressive Strength data set. For details see Section \ref{sec:results}.}
\label{fig:concrete}
\subfigure[SMC $\sim1.2\times10^6$eval'ns of $\pi$]{
\includegraphics[angle=-90,width=\histwidth]{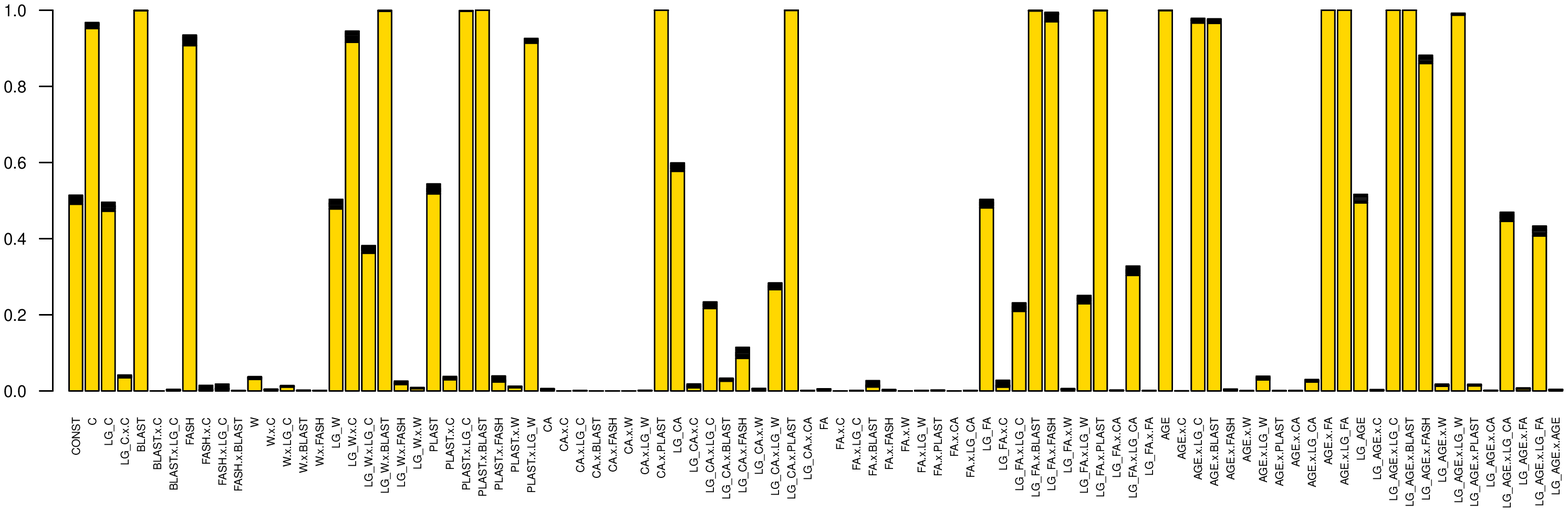}
\label{fig:concrete_smc}
}
\subfigure[AMCMC $2.5\times10^6$eval'ns of $\pi$]{
\includegraphics[angle=-90,width=\histwidth]{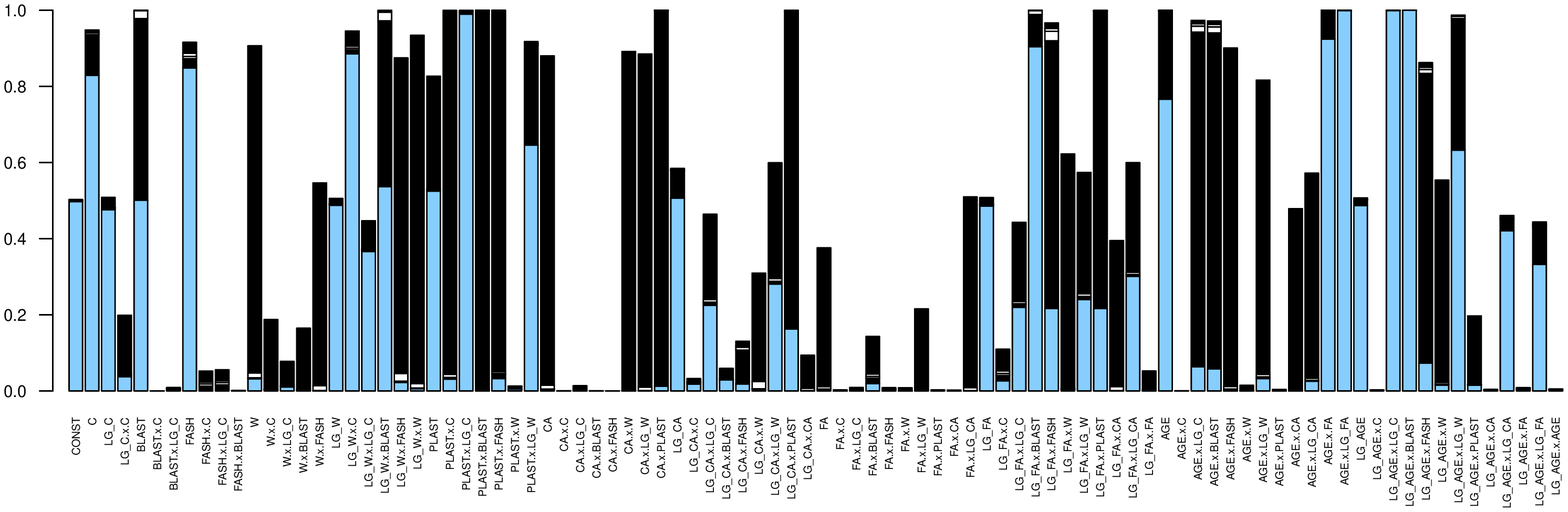}
\label{fig:concrete_amcmc}
}
\subfigure[MCMC $2.5\times10^6$eval'ns of $\pi$]{
\includegraphics[angle=-90,width=\histwidth]{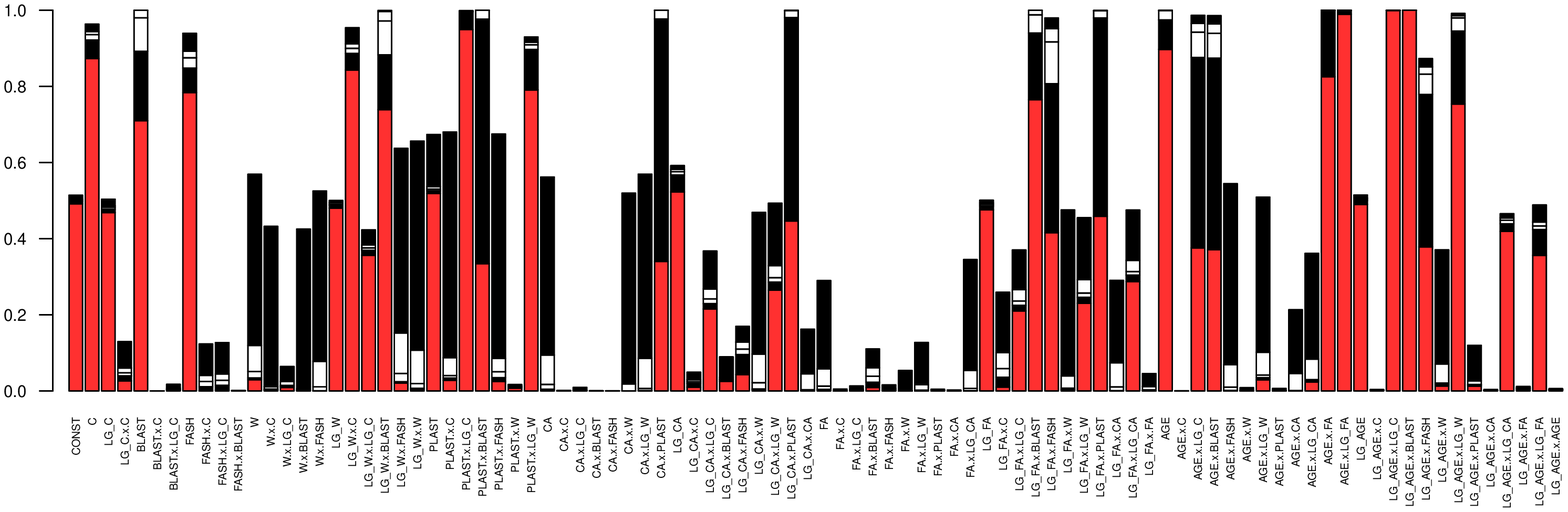}
\label{fig:concrete_mcmc}
}
\vspace{1mm}
\begin{flushleft}
\small{\bfseries Table.} Concrete Compressive Strength data set. Averaged key indicators complementary to Figure \ref{fig:concrete}.
\end{flushleft}
\begin{center}
\begin{tabular}[hb]{l|rrr}
& Sequential MC & Adaptive MCMC & Standard MCMC \\
\noalign{\smallskip}\hline\noalign{\smallskip}
computational time & $0:29:01$ min & $2:02:06$ min & $0:43:17$ min \\
evaluations of $\pi$ & $1.19\times10^6$ & $2.50\times10^6$ & $2.50\times10^6$ \\
average acceptance rate & $30.7\%$ &  $70.4\%$  &  $7.20\%$ \\
length $t$ of the chain $\v x_t$ & & $2.43\times10^7$ & $2.50\times10^6$ \\
moves $\v x_t\neq \v x_{t-1}$ & & $1.76\times10^6$ & $1.79\times10^5$ \\
\end{tabular}
\end{center}
\end{figure*}
\end{center}

\begin{center}
\begin{figure*}
\vspace{-2mm}
\caption{Concrete Compressive Strength data set with main effect restrictions. For details see Section \ref{sec:results}.}
\label{fig:me concrete}
\subfigure[SMC $\sim2.4\times10^6$eval'ns of $\pi$]{
\includegraphics[angle=-90,width=\histwidth]{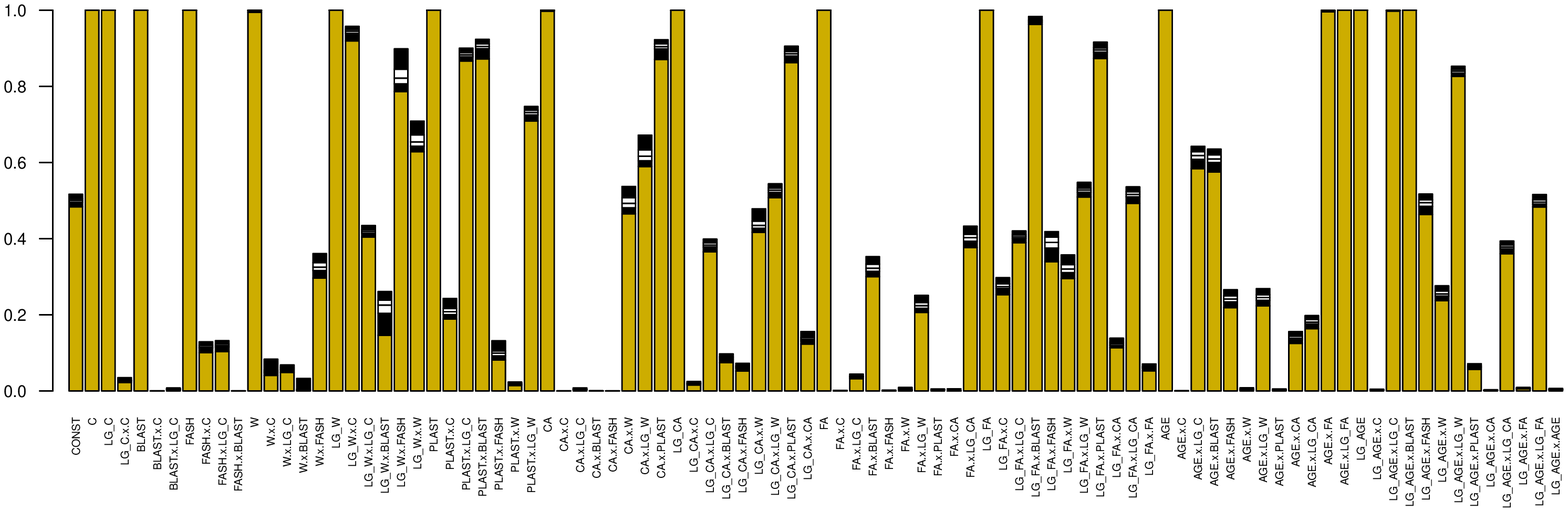}
\label{fig:me concrete_smc}
}
\subfigure[AMCMC $2.5\times10^6$eval'ns of $\pi$]{
\includegraphics[angle=-90,width=\histwidth]{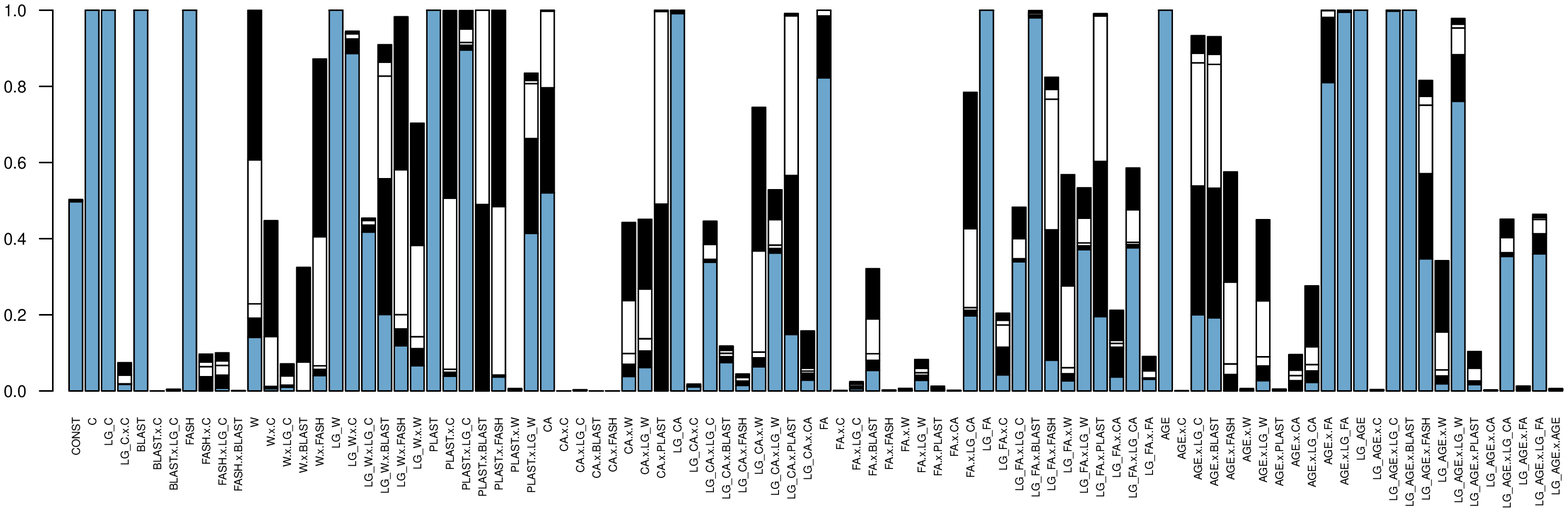}
\label{fig:me concrete_amcmc}
}
\subfigure[MCMC $2.5\times10^6$eval'ns of $\pi$]{
\includegraphics[angle=-90,width=\histwidth]{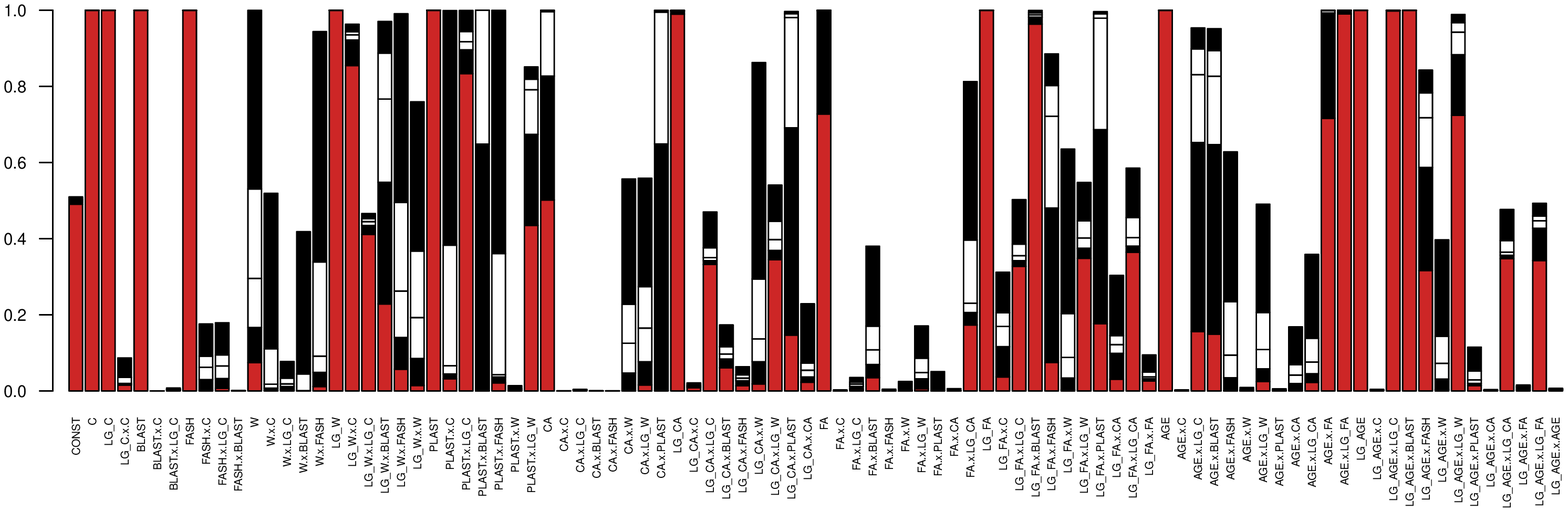}
\label{fig:me concrete_mcmc}
}
\vspace{1mm}
\begin{flushleft}
\small{\bfseries Table.} Concrete Compressive Strength data set with main effect restrictions. Averaged key indicators complementary to Figure \ref{fig:me concrete}.
\end{flushleft}
\begin{center}
\begin{tabular}[hb]{l|rrr}
& Sequential MC & Adaptive MCMC & Standard MCMC \\
\noalign{\smallskip}\hline\noalign{\smallskip}
computational time & $0:43:01$ min & $2:29:16$ min & $0:41:48$ min \\
evaluations of $\pi$ & $2.42\times10^6$ & $2.50\times10^6$ & $2.50\times10^6$ \\
average acceptance rate & $30.98\%$ &  $61.1\%$  &  $5.31\%$ \\
length $t$ of the chain $\v x_t$ & & $2.72\times10^7$ & $2.50\times10^6$ \\
moves $\v x_t\neq \v x_{t-1}$ & & $1.53\times10^6$ & $1.32\times10^5$ \\
\end{tabular}
\end{center}
\end{figure*}
\end{center}

%


%

\begin{center}
\begin{figure*}
\vspace{-2mm}
\caption{Protein data set. For details see Section \ref{sec:results}.}
\label{fig:protein}
\subfigure[SMC $\sim6.1\times10^5$eval'ns of $\pi$]{
\includegraphics[angle=-90,width=\histwidth]{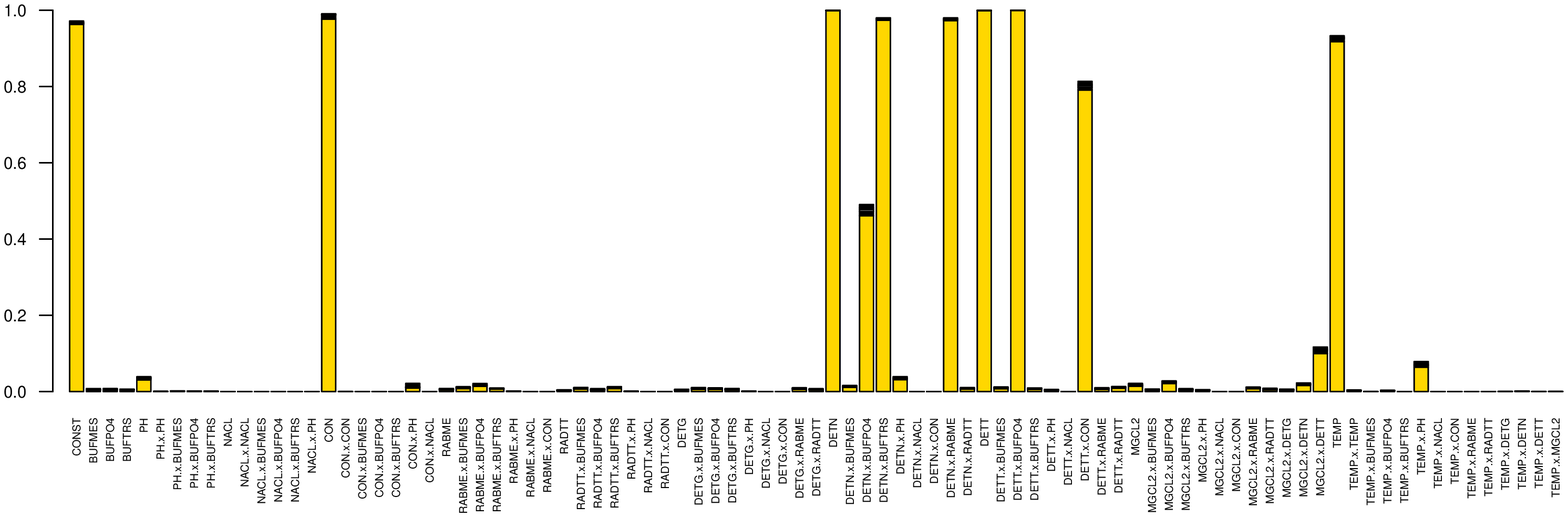}
\label{fig:protein_smc}
}
\subfigure[AMCMC $2.5\times10^6$eval'ns of $\pi$]{
\includegraphics[angle=-90,width=\histwidth]{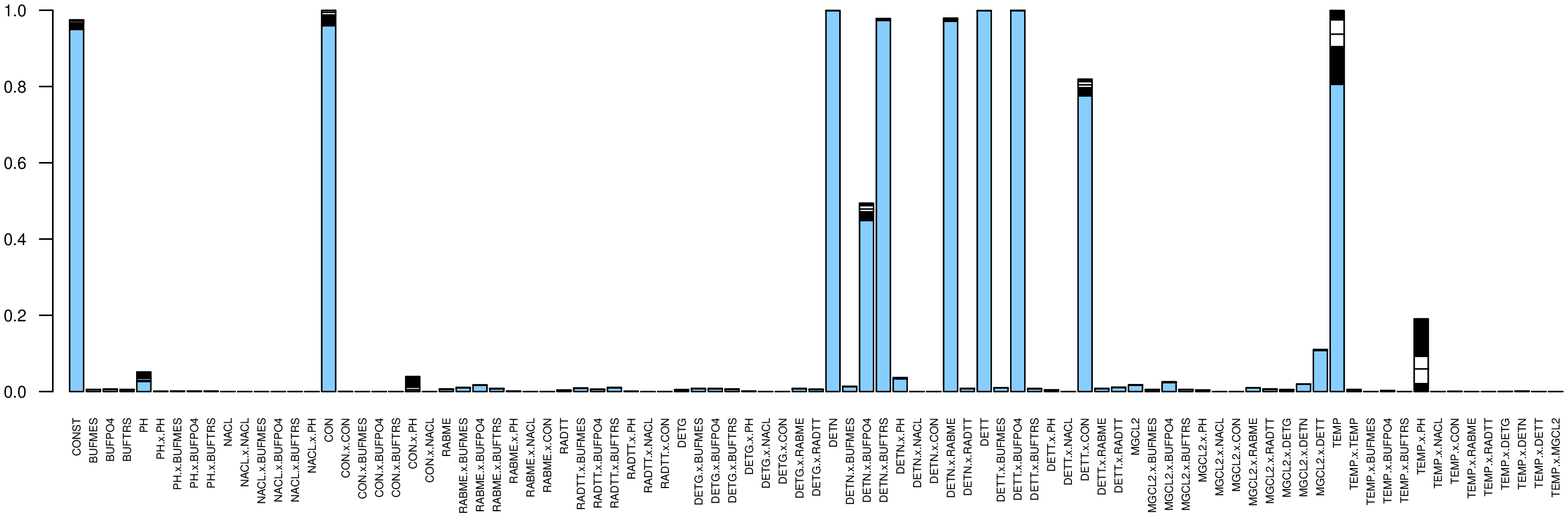}
\label{fig:protein_amcmc}
}
\subfigure[MCMC $2.5\times10^6$eval'ns of $\pi$]{
\includegraphics[angle=-90,width=\histwidth]{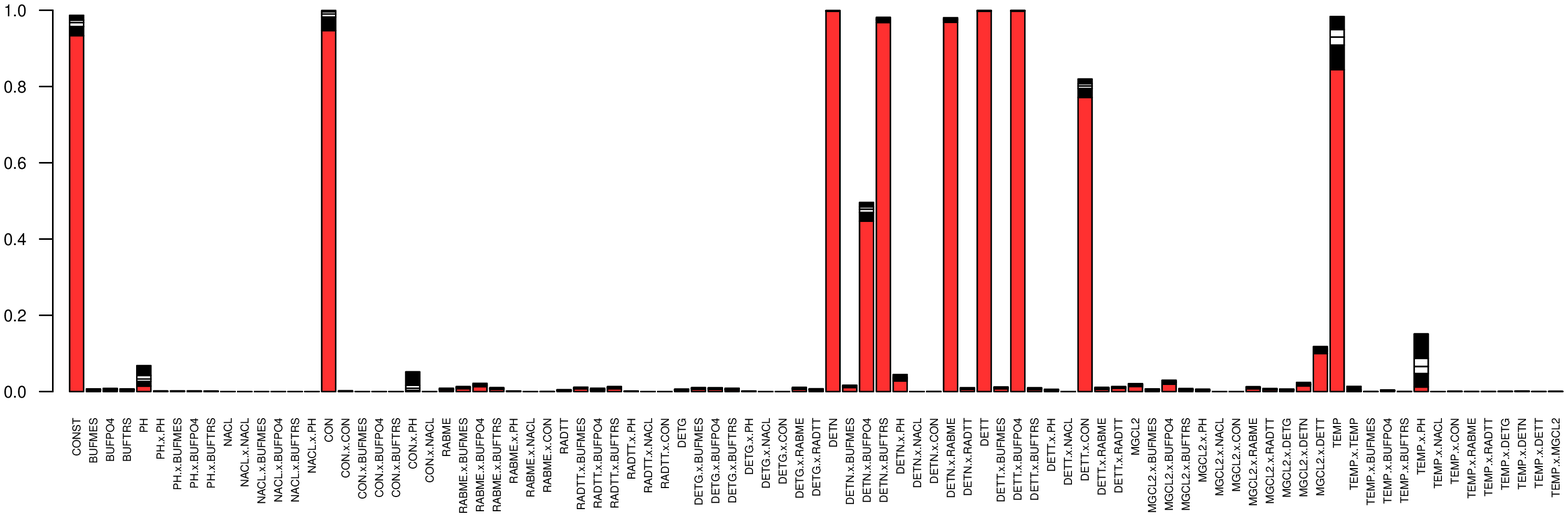}
\label{fig:protein_mcmc}
}
\vspace{1mm}
\begin{flushleft}
\small{\bfseries Table.} Protein data set. Averaged key indicators complementary to Figure \ref{fig:protein}.
\end{flushleft}
\begin{center}
\begin{tabular}[hb]{l|rrr}
& Sequential MC & Adaptive MCMC & Standard MCMC \\
\noalign{\smallskip}\hline\noalign{\smallskip}
computational time & $0:14:55$ min & $3:58:32$ min & $0:29:38$ min \\
evaluations of $\pi$ & $6.17\times10^5$ & $2.50\times10^6$ & $2.50\times10^6$ \\
average acceptance rate & $30.7\%$ &  $60.7\%$  &  $1.20\%$ \\
length $t$ of the chain $\v x_t$ & & $9.19\times10^7$ & $2.50\times10^6$ \\
moves $\v x_t\neq \v x_{t-1}$ & & $1.51\times10^6$ & $3.03\times10^5$ \\
\end{tabular}
\end{center}
\end{figure*}
\end{center}

\begin{center}
\begin{figure*}
\vspace{-2mm}
\caption{Protein data set with main effect restrictions. For details see Section \ref{sec:results}.}
\label{fig:me protein}
\subfigure[SMC $\sim6.1\times10^5$eval'ns of $\pi$]{
\includegraphics[angle=-90,width=\histwidth]{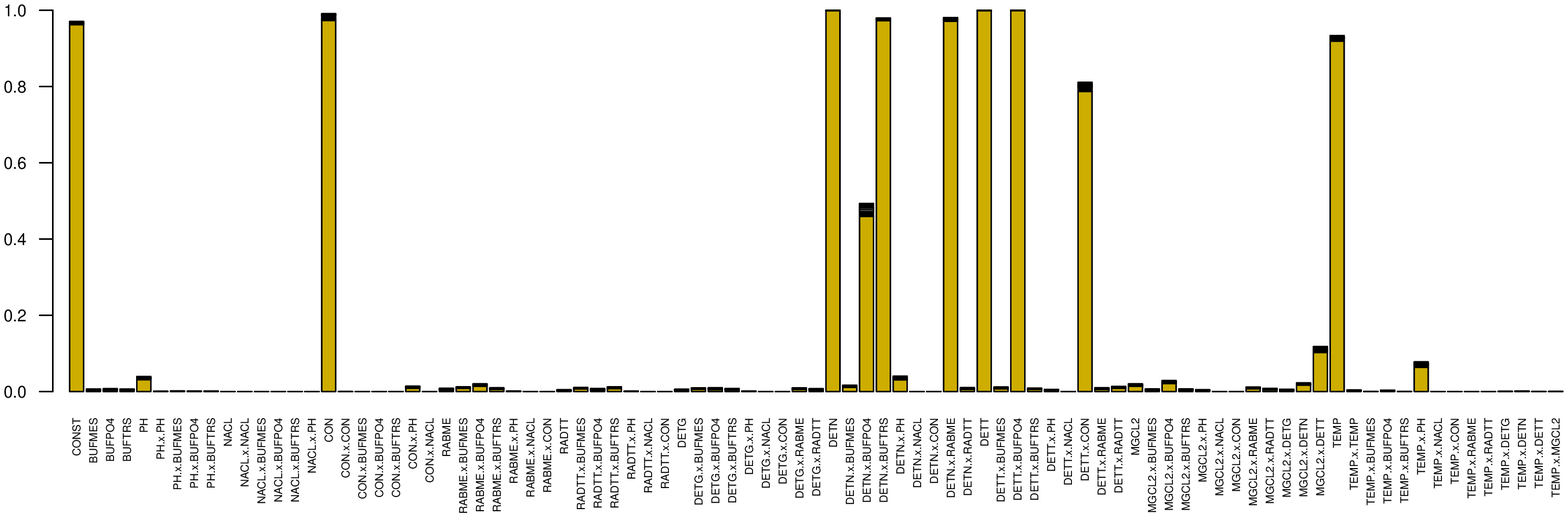}
\label{fig:me protein_smc}
}
\subfigure[AMCMC $2.5\times10^6$eval'ns of $\pi$]{
\includegraphics[angle=-90,width=\histwidth]{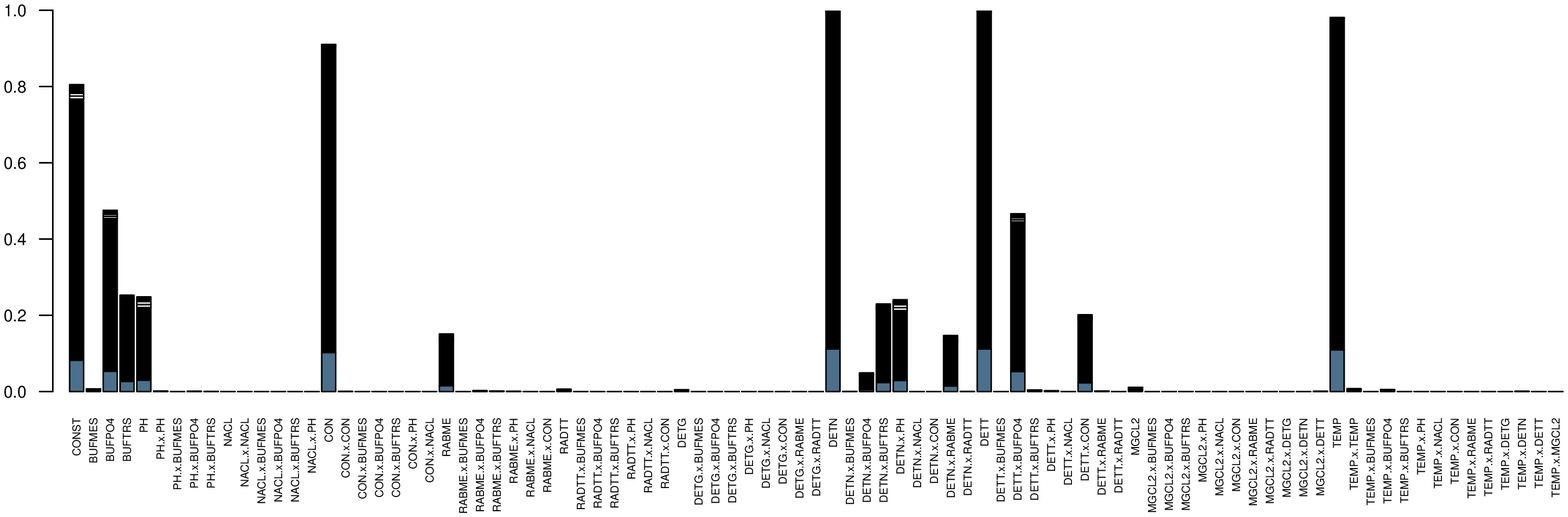}
\label{fig:me protein_amcmc}
}
\subfigure[MCMC $2.5\times10^6$eval'ns of $\pi$]{
\includegraphics[angle=-90,width=\histwidth]{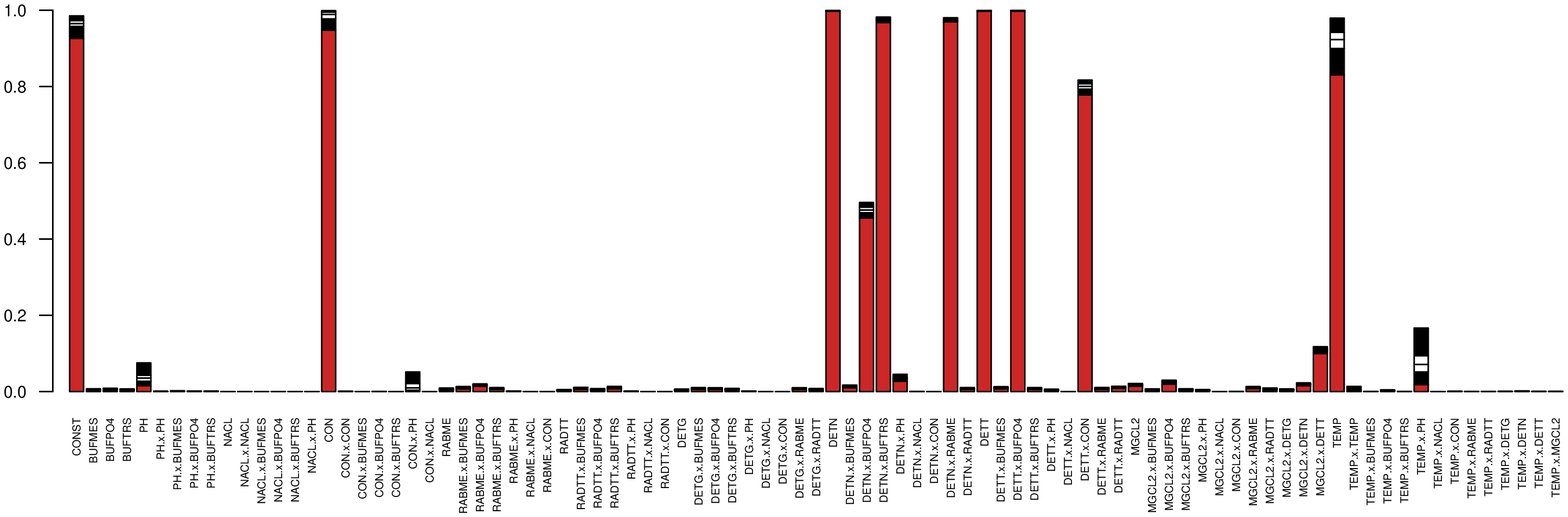}
\label{fig:me protein_mcmc}
}
\vspace{1mm}
\begin{flushleft}
\small{\bfseries Table.} Protein data set with main effect restrictions. Averaged key indicators complementary to Figure \ref{fig:me protein}.
\end{flushleft}
\begin{center}
\begin{tabular}[hb]{l|rrr}
& Sequential MC & Adaptive MCMC & Standard MCMC \\
\noalign{\smallskip}\hline\noalign{\smallskip}
computational time & $0:14:45$ min & $3:32:06$ min & $0:30:21$ min \\
evaluations of $\pi$ & $6.19\times10^5$ & $2.50\times10^6$ & $2.50\times10^6$ \\
average acceptance rate & $26.65\%$ & $22.3\%$ & $1.20\%$ \\
length $t$ of the chain $\v x_t$ & & $1.07\times10^8$ & $2.50\times10^6$ \\
moves $\v x_t\neq \v x_{t-1}$ & & $5.56\times10^6$ & $3.03\times10^5$ \\
\end{tabular}
\end{center}
\end{figure*}
\end{center}

\section*{Acknowledgements}
N.~Chopin is supported by the ANR grant ANR-008-BLAN-0218 ``BigMC'' of the French Ministry of research.

We would like to thank Pierre Jacob and two anonymous referees for their valuable comments on this paper. We acknowledge the StatLib data archive and the UCI Machine Learning Repository for providing the data sets used in this work.


\end{document}